\documentclass[11pt, oneside]{article}

\usepackage[dvips]{graphicx}
\usepackage{epsfig}
\usepackage{tikz}
\usepackage{amsmath}
\usepackage{amssymb}
\usepackage{xcolor}
\usepackage{mathtools}

\usepackage{booktabs}
\usepackage{multicol}
\usepackage{multirow}
\usepackage[labelfont=bf]{caption}
\usepackage{subcaption}

\newcommand{\x}{\mathbf{x}}

\begin{document}
\baselineskip=2pc

\vspace{.5in}

\begin{center}

{\large\bf
Sparse grid implementation of a fixed-point fast sweeping WENO scheme
for Eikonal equations
\footnote{Research was partially supported by NSF grant DMS-1620108.}
}

\end{center}

\vspace{.1in}

\centerline{Zachary M. Miksis\footnote{Department of Applied and Computational Mathematics and Statistics,
University of Notre Dame, Notre Dame, IN 46556, USA. E-mail: zmiksis@nd.edu
},
Yong-Tao Zhang\footnote{Corresponding author. Department of Applied and Computational Mathematics and Statistics,
University of Notre Dame, Notre Dame, IN 46556, USA. E-mail: yzhang10@nd.edu
}
}

\vspace{.3in}

\abstract{Fixed-point fast sweeping methods are a class of explicit iterative methods developed in the literature to efficiently solve steady state solutions of hyperbolic partial differential
equations (PDEs). As other types of fast sweeping schemes, fixed-point fast sweeping methods use the
Gauss-Seidel iterations and alternating sweeping strategy to cover characteristics of hyperbolic PDEs in a certain direction simultaneously in each sweeping order. The resulting iterative schemes have fast convergence
rate to steady state solutions. Moreover, an advantage of fixed-point fast sweeping methods over other types of fast sweeping methods is that they are explicit and do not involve inverse operation of any nonlinear local system. Hence they are robust and flexible, and have been combined with high order accurate weighted essentially non-oscillatory (WENO) schemes to solve various hyperbolic PDEs in the literature. For multidimensional nonlinear problems, high order fixed-point fast sweeping WENO methods still require quite large amount of computational costs. In this technical note, we apply sparse-grid techniques, an effective approximation tool for multidimensional problems, to  fixed-point fast sweeping WENO method for reducing its computational costs. Here we focus on a robust Runge-Kutta (RK) type fixed-point fast sweeping WENO scheme with third order accuracy (Zhang et al. 2006 \cite{ZZC2006}), for solving Eikonal equations, an important class of static Hamilton-Jacobi (H-J) equations.
Numerical experiments on solving multidimensional Eikonal equations and a more general static H-J equation are performed to show that the sparse grid computations of the fixed-point fast sweeping WENO scheme  achieve large savings of CPU times on refined meshes, and at the same time maintain comparable accuracy and resolution with those on corresponding regular single grids.
}

\vspace{.3in}


\pagenumbering{arabic}

\section{Introduction}

In this technical note, we study an efficient approach to reduce the computational costs for solving
the multidimensional Eikonal equations
\begin{equation}
\left\{ \begin{matrix*}[l]
|\nabla \phi(\x)| = f(\x), & & \x \in \Omega \backslash \Gamma \subset \mathbb{R}^d, \\
\phi(\x) = g(\x), & & \x \in \Gamma \subset \Omega,
\end{matrix*}
\right.
\label{eq:eikonal}
\end{equation}
where $\Omega$ is a $d$-dimension computational domain in $R^d$
and $\Gamma$ is a subset of $\Omega$. The given functions $f(\x)$ and $g(\x)$ are Lipschitz continuous, and $f(\x)$ is positive. The Eikonal equations are a very important class of static
Hamilton-Jacobi (H-J) equations \cite{CL}
\begin{equation}
\left\{ \begin{matrix*}[l]
H(\x,\nabla \phi (\x))= f(\x), & & \x \in \Omega \backslash \Gamma \subset \mathbb{R}^d, \\
\phi(\x) = g(\x), & & \x \in \Gamma \subset \Omega,
\end{matrix*}
\right. \label{eq:SHJ}
\end{equation}
where $H$ is the Hamiltonian.
The numerical computations of Eikonal equations
appear in many applications, such as optimal control,
image processing and computer vision, geometric optics,
seismic waves, level set methods, etc.

Due to nonlinearity of the equations and possible singularities in their solutions, it is challenging to design efficient and high order accurate numerical methods
for solving static H-J equations such as the Eikonal equations (\ref{eq:eikonal}). In the literature, a popular approach is to
discretize (\ref{eq:eikonal}) into a nonlinear system and
then design a fast numerical method to solve the
nonlinear system.  Among such methods are the fast marching method
and the fast sweeping method. The fast marching method uses the Dijkstra's algorithm \cite{DI}
and updates the solution by following the Eikonal equations' causality sequentially, e.g., see \cite{Se, SV01, SV03}. In the fast sweeping method \cite{ZOMK, Z, KOQ,
QZZ1, QZZ2, FLZ}, Gauss-Seidel
iterations with alternating orderings are combined with upwind
schemes. Different from the fast marching method, the fast
sweeping method is an iterative method and follows the Eikonal equations' causality
along characteristics in a parallel way, i.e., each Gauss-Seidel iteration with a specific
sweeping ordering covers a family of characteristics in a certain
direction simultaneously.

The iterative framework of fast sweeping method provides certain flexibility to incorporate
high order accuracy schemes for hyperbolic PDEs, such as weighted essentially nonoscillatory (WENO) methods \cite{ZZQ, Xiong} or discontinuous Galerkin (DG) \cite{LSZZ, ZCLZS, WZ} methods, into it for developing high
order fast sweeping methods. In \cite{ZZC2006}, fixed-point fast sweeping WENO methods were designed to solve static H-J equations. Different from other fast sweeping methods, fixed-point fast sweeping methods adopt the Gauss-Seidel idea and alternating sweeping strategy to the time-marching type fixed-point iterations. They are explicit schemes and do not involve inverse operation of nonlinear local systems which have to be done in other types of fast sweeping methods, hence are much more easier to be applied in solving various hyperbolic equations using any monotone numerical fluxes and high order nonlinear WENO approximations. For example, how to efficiently solve steady state problems of hyperbolic conservation laws is important and challenging \cite{ChouS,W.Chen}.  In \cite{WuLiang, LZZ}, fixed-point fast sweeping WENO methods were applied in solving nonlinear hyperbolic conservation laws. Numerical experiments performed in \cite{ZZC2006, WuLiang, LZZ} show that more than $50\%$ computational costs are saved by using fixed-point fast sweeping methods rather than direct time-marching methods to converge to steady states of high order WENO schemes.

Since high order WENO methods require more operations than
many other schemes due to their sophisticated nonlinearity and high order accuracy, the associated computational costs increase significantly when the number of grid points is large for multidimensional problems.
Sparse-grid techniques, an efficient approach for solving high-dimensional problems, have been developed in the literature to reduce the number of grid points needed in the simulations. See \cite{BG, JGSG} for a review.
In 1991, sparse-grid techniques were introduced in \cite{Zg} to reduce the number of degrees of freedom in finite element method. As an approach for practical implementation of sparse-grid techniques, the sparse-grid combination technique was developed in \cite{GSZ}. The main idea of sparse-grid combination technique is to compute the final solution as a linear combination of solutions on semi-coarsened grids, and the coefficients of the linear combination are taken such that there is a canceling in leading-order error terms and the resulting accuracy order is kept to be the same as that on a single full grid. The sparse-grid combination technique was applied to linear schemes in \cite{LKV1,LKV2} in early time. Recently it has been applied to nonlinear WENO schemes in \cite{LuZhang1, LCZ, ZZ21} for solving hyperbolic conservation laws or convection-diffusion equations, where numerical results show that significant computational times are saved, while both accuracy and stability of the nonlinear WENO schemes are maintained for simulations on sparse grids. In this technical note, we follow the way in our previous work and apply the sparse-grid combination technique to a fixed-point fast sweeping WENO method for solving multidimensional Eikonal equations. A Runge-Kutta (RK) type fixed-point fast sweeping WENO scheme with third order accuracy developed in \cite{ZZC2006} is used in this paper, since our numerical experiments find that this RK fixed-point fast sweeping WENO scheme is very robust for simulations performed on sparse grids. The rest of the paper is organized as following. In Section 2, we describe the algorithm how to apply the sparse-grid combination technique to the RK fixed-point fast sweeping WENO scheme. In Section 3, various numerical experiments including solving multidimensional Eikonal equations and a more general static H-J equation with smooth or non-smooth solutions, are carried out
to show that the sparse grid computations of the fixed-point fast sweeping WENO scheme save large amount of CPU times, especially on refined meshes, and at the same time maintain comparable simulation results with those on corresponding regular single grids.
Conclusions are given in Section 4.

\section{Description of the numerical algorithm}

In this section, we first review the RK fixed-point fast sweeping WENO scheme in \cite{ZZC2006}, then describe
the algorithm to implement it on sparse grids.

\subsection{The RK fixed-point fast sweeping WENO scheme}

The fixed-point fast sweeping WENO schemes in \cite{ZZC2006} were developed by applying the Gauss-Seidel idea and alternating sweeping strategy to the time-marching schemes to solve the static H-J equations \eqref{eq:SHJ}.
The RK fixed-point fast sweeping WENO scheme used in this paper is based on the second order total variation diminishing (TVD) RK time-marching scheme \cite{Shu} and third order WENO approximations to spatial derivatives.
Here we take the two dimensional case as an example to describe the method, which is similar for higher dimensional cases. The computational domain $\Omega$ is partitioned by a Cartesian grid $\{(x_i, y_j), 1\leq i \leq I, 1\leq j \leq J\}$, with uniform grid sizes $h_x$ and $h_y$ in the $x$ and $y$ directions respectively. Denote the viscosity
numerical solution of \eqref{eq:SHJ} at a grid point $(x_i, y_j)$ by $\phi_{i,j}$.
The RK fixed-point fast sweeping scheme in \cite{ZZC2006} has the following form:
\begin{align}
\phi^{(1)}_{i,j} &= \phi^n_{i,j} + \gamma \left( \frac{1}{ \frac{\alpha_x}{h_x} + \frac{\alpha_y}{h_y} } \right) \left[ f_{i,j}  - \hat{H} ((\phi_x)^-_{i,j},(\phi_x)^+_{i,j}; (\phi_y)^-_{i,j},(\phi_y)^+_{i,j}) \right]; \label{eq:RKs1} \\
\phi^{n+1}_{i,j} &= \phi^{(1)}_{i,j} + \frac{1}{2}\gamma \left( \frac{1}{ \frac{\alpha_x}{h_x} + \frac{\alpha_y}{h_y} } \right) \left[ f_{i,j}  - \hat{H} ((\phi_x)^-_{i,j},(\phi_x)^+_{i,j}; (\phi_y)^-_{i,j},(\phi_y)^+_{i,j}) \right]. \label{eq:RKs2}
\end{align}
Here $\phi^n_{i,j} $ and $\phi^{n+1}_{i,j}$ are the numerical solution values at iteration step $n$ and $n+1$, respectively. $f_{i,j}$ denotes the value of $f$ at a grid point $(x_i, y_j)$.
$\hat{H}$ is a monotone numerical Hamiltonian \cite{OSherShu}.  $(\phi_x)^-_{i,j}$ is an approximation of $\phi_x$
at the grid point $(x_i, y_j)$ when the wind ``blows''
from the left to the right, and $(\phi_x)^+_{i,j}$ is an approximation of $\phi_x$
at the grid point $(x_i, y_j)$ when the wind ``blows'' from the right to the left. It is similar for
$y$-direction approximations
$(\phi_y)^-_{i,j}$ and $(\phi_y)^+_{i,j}$. $\gamma$ is a parameter. To guarantee that the fixed-point iteration is a contractive mapping and converges, suitable values of $\gamma$ need to be taken. In the context of time-marching schemes, $\gamma$ is actually the Courant-Friedrichs-Lewy (CFL) number.
\begin{equation}
\alpha_x  = \max_{\substack{A \leq u \leq B \\ C \leq v \leq D}} |H_1(u,v)|, \qquad \alpha_y  = \max_{\substack{A \leq u \leq B \\ C \leq v \leq D}} |H_2(u,v)|.
\label{alpha}
\end{equation}
$H_i(u,v)$ is the partial derivative of $H$ with respect to the $i$th argument, or  the Lipschitz constant of $H$ with respect to the $i$th argument. $[A,B]$ is the value range for $\phi_x^\pm$, and $[C,D]$ is the value range for $\phi_y^\pm$. For the Eikonal equation (\ref{eq:eikonal}), we have $\alpha_x = \alpha_y = 1$.

\begin{figure}
\centering
\includegraphics[width=4in]{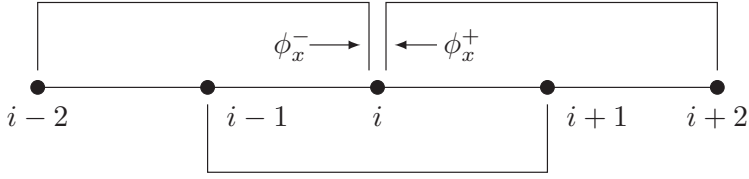}
\caption{Stencils of the third-order WENO approximations for derivatives.}
\label{stencil}
\end{figure}

For first order scheme, simple first order upwind finite difference approximations for $\phi_x$ and $\phi_y$ are used. To obtain a high order scheme, in \cite{ZZC2006} $(\phi_x)^-_{i,j}$, $(\phi_x)^+_{i,j}$, $(\phi_y)^-_{i,j}$, and $(\phi_y)^+_{i,j}$ are computed by a third order WENO scheme, which is also used in \cite{ZZQ}. See Fig. \ref{stencil} for an illustration of the interpolation stencils used. The WENO approximation of $\phi_x$ at the grid point $(x_i,y_j)$ when the wind ``blows'' left-to-right is
\begin{equation}
(\phi_x)^-_{i,j} = (1-w_-) \left( \frac{\phi_{i+1,j} - \phi_{i-1,j}}{2h_x} \right) + w_- \left( \frac{3\phi_{i,j} - 4\phi_{i-1,j} + \phi_{i-2,j}}{2h_x} \right),
\label{weno1}
\end{equation}
where
\begin{equation}
w_- = \frac{1}{1+2r_-^2}, \quad r_- = \frac{\epsilon + (\phi_{i,j} - 2\phi_{i-1,j} + \phi_{i-2,j})^2}{\epsilon + (\phi_{i+1,j} - 2\phi_{i,j} + \phi_{i-1,j})^2};
\end{equation}
when the wind ``blows'' right-to-left, the WENO approximation is
\begin{equation}
(\phi_x)^+_{i,j} = (1-w_+) \left( \frac{\phi_{i+1,j} - \phi_{i-1,j}}{2h_x} \right) + w_+ \left( \frac{-\phi_{i+2,j} + 4\phi_{i+1,j} - 3\phi_{i,j}}{2h_x} \right),
\label{weno2}
\end{equation}
where
\begin{equation}
w_+ = \frac{1}{1+2r_+^2}, \quad r_+ = \frac{\epsilon + (\phi_{i+2,j} - 2\phi_{i+1,j} + \phi_{i,j})^2}{\epsilon + (\phi_{i+1,j} - 2\phi_{i,j} + \phi_{i-1,j})^2}.
\end{equation}
Here $\epsilon$ is a small value to avoid that the denominator becomes zero. The WENO approximations of $\phi_y$ are computed similarly. If we take $w_-=w_+=1/3$ in (\ref{weno1}) and (\ref{weno2}), then third order linear upwind approximations are obtained. In this paper, we use the Lax-Friedrichs numerical Hamiltonian \cite{OSherShu}, which has the following form for a Hamiltonian $H(u,v)$:
\begin{equation}
\hat{H}^{LF}(u^-,u^+;v^-,v^+) = H \left( \frac{u^- + u^+}{2}, \frac{v^- + v^+}{2} \right) - \frac{1}{2} \alpha_x (u^+ - u^-) - \frac{1}{2} \alpha_y (v^+ - v^-),
\end{equation}
where $\alpha_x$ and $\alpha_y$ have the same definition as (\ref{alpha}). Note that the philosophy
of Gauss-Seidel iterations is adopted to compute the approximations for derivatives, namely, we always use the newest available values of $\phi$ in the interpolation stencils
to compute the approximations for $(\phi_x)^-_{i,j}$, $(\phi_x)^+_{i,j}$, $(\phi_y)^-_{i,j}$, $(\phi_y)^+_{i,j}$
in (\ref{eq:RKs1}) and (\ref{eq:RKs2}).

We summarize the RK fixed-point fast sweeping WENO (RK FPFS-WENO) algorithm as the following.
\begin{enumerate}
\item {\em Initialization:} according to the boundary condition $\phi(x,y) = g(x,y)$, $(x,y) \in \Gamma$, assign exact values or interpolated values at grid points whose distances to $\Gamma$ are less than or equal to $(m-1)$ grid sizes, where $m$ is the number of grid points in small stencils of WENO approximations. For example, $m=3$ for the third order WENO approximations used here. These values are fixed during iterations. For robust simulations, the solution from the non-fully-converged (i.e., using a much larger convergence threshold value $\delta$ than that of the WENO sweeping used in the step 3 below; specific values given in the numerical example section) first order sweeping computation (i.e., using the first order upwind approximations for these derivatives in $\hat{H}$ of the scheme \eqref{eq:RKs1}-\eqref{eq:RKs2}) is used as the initial guess at all other grid points, while a big value (e.g. $10$ in this paper) is used as the initial guess for the first order sweeping computation.
\item {\em Iterations:} perform the Gauss-Seidel iterations \eqref{eq:RKs1}-\eqref{eq:RKs2} with four alternating direction sweepings:
\begin{align*}
&(a) \, i=1:I, \, j=1:J; \\
&(b) \, i = I:1, \, j = 1:J; \\
&(c) \, i=I:1, \, j = J:1; \\
&(d) \, i = 1:I, \, j = J:1.
\end{align*}
Each sweeping direction is completed in full for the first  Runge-Kutta stage before moving to the second Runge-Kutta stage, and the sweeping direction
of both stages should be same during one sweeping. High order extrapolations are used for the ghost points when calculating the high order WENO approximations of the derivatives for grid points on the boundary of the computational domain, as in \cite{ZZQ}.

\item {\em Convergence:} if
$$
\| \phi^{n+1} - \phi^{n} \|_{L^\infty} \leq \delta,
$$
where $\delta$ is a given convergence threshold value and $\| \cdot \|_{L^\infty}$ denotes the $L^\infty$ norm, the algorithm converges and we stop the iterations.
\end{enumerate}


\subsection{RK FPFS-WENO scheme on sparse grids}

\begin{figure}
\centering
\includegraphics[width=4.5in]{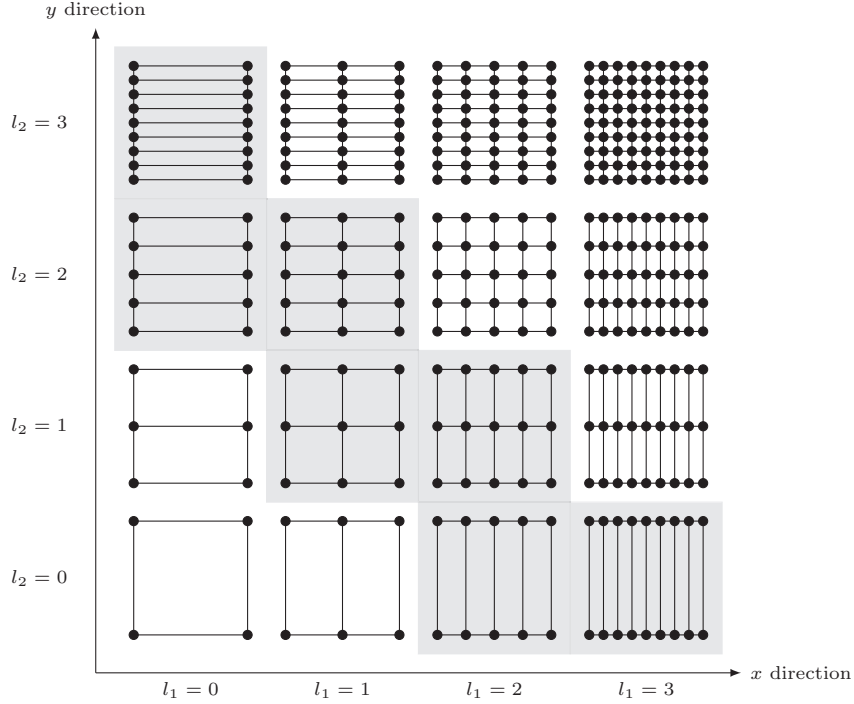}
\caption{Illustration of 2D sparse grids $\{\Omega^{l_1,l_2}\}$ for one cell of a root grid. Here the cell indicated by the levels $l_1=0, l_2=0$ is one cell of the root grid $\Omega^{0,0}$, and the side length of the cell is $H$. The finest level $N_L=3$. Highlighted grids are those on which PDEs are solved.}
\label{sparse_grids}
\end{figure}

In this section, we describe how to implement the RK FPFS-WENO method on sparse grids by using the sparse-grid combination technique, for improving the method's efficiency in solving multidimensional Eikonal equations. Here two dimensional (2D) cases are used to illustrate the idea. Algorithm procedures for higher dimensional cases are similar. We consider a square computational domain $[a,b]^2$ for simplicity of the description, and construct semi-coarsened sparse grids as the following. Note that the procedure here can be applied to any rectangular domain straightforwardly.
The domain is first partitioned into the coarsest  grid $\Omega^{0,0}$ with $N_r$ cells in each direction and mesh size $H = \frac{b-a}{N_r}$. $\Omega^{0,0}$ is called a root grid. Then a multi-level refinement on the root grid is done to construct a family of semi-coarsened grids $\{ \Omega^{l_1,l_2} \}$, with mesh sizes $h_{l_1} = 2^{-l_1}H$ in the $x$-direction and $h_{l_2} = 2^{-l_2}H$ in the $y$-direction, where $l_1 = 0,1,\dots,N_L$ and $l_2 = 0,1,\dots,N_L$.  The superscripts $l_1, l_2$ are the refinement levels relative to the root grid $\Omega^{0,0}$ in the $x$ and the $y$ directions respectively, and $N_L$ is the finest level. Here the finest grid is $\Omega^{N_L,N_L}$ with the mesh size $h=2^{-N_L}H$ in both $x$ and $y$ directions. Actually $\Omega^{N_L,N_L}$ is corresponding to a single full grid in regular single-grid computations. Figure \ref{sparse_grids} is an
illustration of 2D sparse grids $\{\Omega^{l_1,l_2}\}$ for one cell of a root grid, with $N_L=3$.
We apply the spare-grid combination techniques. The Eikonal equation (\ref{eq:eikonal}) is not directly solved by the RK FPFS-WENO method on a single full grid $\Omega^{N_L,N_L}$, but on the set $\{\Omega^{l_1,l_2}\}_I$ of the following  $(2N_L + 1)$ sparse grids:
$$
\left\{ \Omega^{0,N_L}, \Omega^{1,N_L - 1}, \cdots, \Omega^{N_L - 1,1}, \Omega^{N_L,0}  \right\} \quad \text{and} \quad \left\{ \Omega^{0,N_L-1}, \Omega^{1,N_L - 2}, \cdots, \Omega^{N_L - 2,1}, \Omega^{N_L-1,0}  \right\},
$$
with $I$ being the index set
$$
I = \{(l_1,l_2)|l_1 + l_2 = N_L \quad \text{or} \quad l_1+l_2 = N_L - 1\}.
$$
Then we have $(2N_L+1)$ sets of numerical solutions $\{\Phi^{l_1,l_2}\}_I$, where each set of numerical solutions is corresponding to each sparse grid in $\{\Omega^{l_1,l_2}\}_I$. The final step is to combine these sparse grid solutions $\{\Phi^{l_1,l_2}\}_I$ to obtain the final solution on the finest grid $\Omega^{N_L,N_L}$.
This is implemented by first using a prolongation operator $P^{N_L,N_L}$ to map each sparse grid solution $\Phi^{l_1,l_2}$ onto the finest grid $\Omega^{N_L,N_L}$, then combining these solutions to form the final solution $\hat{\Phi}^{N_L,N_L}$ on $\Omega^{N_L,N_L}$. Next we describe the prolongation technique in details, then summarize the whole algorithm.

\subsubsection{Prolongation operator and WENO interpolation}

Given the numerical solution $\Phi^{l_1,l_2}$ on $\Omega^{l_1,l_2}$, a prolongation operator $P^{N_L,N_L}$
generates numerical values $P^{N_L,N_L}\Phi^{l_1,l_2}$ for all grid points on $\Omega^{N_L,N_L}$.
Prolongation is usually implemented by interpolation procedure. Studies in \cite{GSZ,LKV1,LKV2} for linear schemes and in \cite{LuZhang1,LCZ,ZZ21} for nonlinear schemes show that the final solution resulted from the spare-grid combination techniques can achieve the similar accuracy orders as the numerical schemes, as long as the accuracy order of interpolations in the prolongations is not less than the accuracy order of the numerical schemes used to solve PDEs on sparse grids. Hence we use third order interpolations here for prolongations.
If solutions are smooth, simple Lagrange interpolation can be used directly. The interpolations are carried out in the dimension by dimension way. In a two dimensional domain, first $(N_r2^{l_1-1})$ quadratic polynomials $P_i^2(x)$, $i = 1, \cdots, N_r2^{l_1-1}$, are constructed along the $x$-direction grid lines using third order Lagrange interpolation. Three adjacent grid points are used in each interpolation. Each polynomial $P_i^2(x)$ is then evaluated on the grid points of $\Omega^{N_L,l_2}$ (the most refined grid in the $x$-direction). Then the same interpolation procedure is performed in every grid line of the $y$ direction with a fixed $x$-coordinate on the grid $\Omega^{N_L,l_2}$, and the obtained polynomials are evaluated on the grid points of $\Omega^{N_L,N_L}$ to get
$P^{N_L,N_L}\Phi^{l_1,l_2}$.


Because solutions of H-J equations may develop  discontinuous derivatives and not be smooth, it is more robust to use WENO interpolations in the prolongation for a general case. Here a third order WENO interpolation is used and detailed formulas are given as following.  We describe the interpolation for a $x$-direction grid line, and it is similar for $y$-direction.  Given numerical values $\phi_{i-1,j}$, $\phi_{i,j}$ and $\phi_{i+1,j}$ at the
grid points $x_{i-1}$, $x_i$ and $x_{i+1}$ along the line $y=y_j$, we compute the third order WENO interpolation
$\phi_{WENO}(x)$ for any point $x \in [x_{i-1/2}, x_{i+1/2})$, where $x_{i-1/2}=(x_{i-1}+x_i)/2$ and $x_{i+1/2}=(x_{i}+x_{i+1})/2$. Let $h$ be the grid size of the uniform mesh, we write the point $x$ as $x = x_{i-1} + \tilde{\alpha}h$ with $\tilde{\alpha} \in [1/2,3/2)$. The WENO interpolation is
\begin{equation}
\phi_{WENO}(x) = w_1 P^1_{(1)}(x) + w_2 P^1_{(2)}(x),
\end{equation}
where $P^1_{(1)}(x)$ and $P^1_{(2)}(x)$ are second order approximations computed as
\begin{equation}
P^1_{(1)}(x) = \tilde{\alpha}\phi_{i,j} - (\tilde{\alpha} - 1)\phi_{i-1,j}, \qquad P^1_{(2)}(x) = (\tilde{\alpha} - 1)\phi_{i+1,j} - (\tilde{\alpha} - 2)\phi_{i,j}.
\end{equation}
The nonlinear weights $w_1$ and $w_2$ are computed as
\begin{equation}
w_1 = \frac{\tilde{w}_1}{\tilde{w}_1 + \tilde{w}_2}, \qquad w_2 = 1 - w_1,
\end{equation}
with
\begin{equation}
\tilde{w}_1 = \frac{\gamma_1}{(\epsilon + \beta_1)^2}, \qquad \tilde{w}_2 = \frac{\gamma_2}{(\epsilon + \beta_2)^2},
\end{equation}
where $\gamma_1 = 1 - \tilde{\alpha}/2$, $\gamma_2 = \tilde{\alpha}/2$, $\beta_1 = (\phi_{i,j} - \phi_{i-1,j})^2$, and $\beta_2 = (\phi_{i+1,j} - \phi_{i,j})^2$.
$\epsilon$ is a small
positive number used to avoid the denominator becoming $0$, and its value is specified in the next numerical
experiment section.

\subsubsection{Algorithm summary}

We summarize the algorithm of the RK FPFS-WENO scheme on sparse grids as following.

\bigskip
\noindent\textbf{Algorithm: sparse grid RK FPFS-WENO scheme}
\begin{enumerate}
\item {\em Restriction step:} perform the initialization step of the RK FPFS-WENO algorithm in section 2.1 on the aforementioned  $(2N_L + 1)$ sparse grids $\{\Omega^{l_1,l_2}\}_I$.

\item {\em Sweeping step:} on each sparse grid $\Omega^{l_1,l_2}$ in $\{\Omega^{l_1,l_2}\}_I$, perform the RK FPFS-WENO iterations to solve the Eikonal equation (\ref{eq:eikonal}). Then we produce $(2N_L +1)$ sets of converged solutions $\{\Phi^{l_1,l_2}\}_I$ for the Eikonal equation (\ref{eq:eikonal}).

\item {\em Prolongation step:} on each sparse grid $\Omega^{l_1,l_2}$ in $\{\Omega^{l_1,l_2}\}_I$, use the prolongation operator $P^{N_L,N_L}$ on $\Phi^{l_1,l_2}$ to map it onto the most refined grid $\Omega^{N_L,N_L}$, and obtain the solution $P^{N_L,N_L}\Phi^{l_1,l_2}$.

\item {\em Combination step:} compute the final solution $\hat{\Phi}^{N_L,N_L}$ by taking the combination
\begin{equation}
\hat{\Phi}^{N_L,N_L} = \sum_{l_1 + l_2 = N_L} P^{N_L,N_L} \Phi^{l_1.l_2} - \sum_{l_1+l_2 = N_L-1} P^{N_L,N_L}\Phi^{l_1,l_2}.
\end{equation}
\end{enumerate}

In three dimensional (3D) or higher dimensional cases, the algorithm follows similar procedure while prolongation operations are carried out in additional spatial directions. The sparse-grid combination formula for higher dimensional problems is provided in the literature, e.g. \cite{GSZ}. In this technical notes, the following 3D formula is also used:
\begin{align}
\hat{\Phi}^{N_L,N_L,N_L} &= \sum_{l_1 + l_2 + l_3= N_L} P^{N_L,N_L,N_L} \Phi^{l_1.l_2,l_3} - 2\sum_{l_1+l_2 +l_3= N_L-1} P^{N_L,N_L,N_L}\Phi^{l_1,l_2,l_3} \nonumber \\
&+ \sum_{l_1 + l_2 + l_3= N_L-2} P^{N_L,N_L,N_L} \Phi^{l_1.l_2,l_3}.
\end{align}




\section{Numerical Examples}

In this section, we perform numerical experiments on solving multidimensional Eikonal equations to test the sparse grid RK FPFS-WENO method and show a large amount of CPU time savings by comparisons with corresponding single-grid simulations. Although theoretical error analysis on linear schemes for linear PDEs \cite{GSZ, LKV1} has been carried out to show that the sparse-grid combination leads to a canceling in leading-order errors of numerical solutions on semi-coarsened sparse grids, hence the accuracy order of the final solution of a sparse-grid computation is kept to be almost the same as that on the corresponding single-grid simulation, such sparse grid error analysis is very difficult to carry out for the WENO methods due to their high nonlinearity. Following our previous studies \cite{LuZhang1, LCZ, ZZ21}, numerical experiments are used to verify the third order accuracy for the sparse grid RK FPFS-WENO scheme in this note, rather than theoretical analysis. Specifically, mesh refinement studies are carried out to compute numerical convergence rates on successively refined grids, for problems with smooth solutions. In \cite{ZZ21}, two different approaches, ``refine root grid'' and ``refine levels'', are studied for mesh refinement in sparse-grid computations. For example, for 3D sparse grids with a $10 \times 10\times 10$ root grid and
$N_L=3$, the finest grid is $80 \times 80\times 80$. The ``refine root grid'' approach is to refine the root grid, while the total number of semi-coarsened sparse-grid levels $N_L+1$ is kept unchanged. So if the root grid is refined once to be $20 \times 20\times 20$, we obtain the finest grid $160 \times 160\times 160$. The ``refine levels'' approach refines the sparse-grid levels, while keeping the root grid fixed. So if $N_L=3$ is refined once to be $N_L=4$, with the fixed $10 \times 10\times 10$ root grid, the finest grid $160 \times 160\times 160$ is also obtained. It is discovered in \cite{ZZ21} that although the ``refine levels'' approach is more efficient and saves more CPU time costs than the ``refine root grid'' approach, it has obvious accuracy order reductions for the nonlinear sparse grid WENO schemes. The ``refine root grid'' approach can always achieve the desired accuracy order of the sparse grid WENO scheme. Hence in this technical note, we use the ``refine root grid'' approach in mesh refinement studies. $N_L = 3$ is used for all sparse grid computations.

We first test the sparse grid RK FPFS-WENO method on problems with smooth solutions to study its numerical accuracy orders. Then the method is applied to problems with non-smooth solutions to show its nonlinear stability. For all numerical examples, we take $\epsilon = 10^{-6}$ in the WENO scheme for both the iterations and the WENO interpolation in the prolongation operator. The convergence threshold value is taken as $\delta = 10^{-11}$ for the third order WENO sweeping, and we take $\delta = 10^{-4}$ in the non-fully-converged first order sweeping to provide initial values for the WENO sweeping. As in \cite{ZZC2006}, we select the largest $\gamma$ value for each problem that provides the iteration convergence with the fastest speed on all semi-coarsened sparse grids in the sparse-grid combination, for the purpose of testing the computational efficiency of the algorithm. To identify the largest possible $\gamma$ value for a problem, we gradually increase / decrease the value of $\gamma$ from an initial value.
In this section, we use $N_h$ to
denote the number of computational cells in one spatial direction of the most refined grid in sparse grids or the corresponding single grid.


\paragraph{Example 1 (A linear problem with smooth solution).}
Consider the following 2D linear problem
\begin{equation}
\phi_x + \phi_y = 0, \qquad (x,y) \in \Omega \backslash \Gamma,
\end{equation}
where $\Omega = [0,2\pi]^2$ and $\Gamma = \{(x,y) \in \Omega \, | \, x=0 \,\, \text{or} \,\, y = 0\}$.
The inflow boundary conditions are applied on $\Gamma$:
\begin{equation}
\phi(x,0) = \sin(x),  \qquad \phi(0,y) = -\sin(y).
\end{equation}
 This problem has the exact solution
\begin{equation}
\phi(x,y) = \sin(x-y).
\end{equation}
For this linear problem with a smooth solution, we solve it by the sparse grid RK FPFS scheme with $\gamma = 1$ and the third order linear upwind approximations to the derivatives, to verify the error analysis results for linear schemes applied to linear PDEs in the literature e.g. \cite{GSZ,LKV1,LKV2}. The third order Lagrange interpolation for prolongation is employed in sparse-grid computations. We perform simulations on both sparse grids and the corresponding single grids, and compare their results. The $L^1$ errors, $L^\infty$ errors and their numerical accuracy orders and CPU times are reported in
Table~\ref{tab:smooththird}. The third order accuracy is obtained for both sparse-grid computations and the corresponding single-grid ones, along with the mesh refinement.
This is consistent with the error analysis results for linear schemes in solving linear PDEs in \cite{GSZ,LKV1,LKV2}.
Comparing the numerical errors of sparse-grid computations and the corresponding single-grid ones, we observe that
their $L^1$ errors are comparable. The $L^\infty$ errors of sparse-grid computations are larger than the corresponding single-grid computations. In terms of computational efficiency, on refined meshes we see around $50\% \sim 80\%$ CPU time saved for simulations on sparse grids vs single grids, for this example.

\begin{table}
\centering
\begin{tabular}{ccccccc}
\toprule[1.5pt]
\multicolumn{7}{c}{Single grid} \\
\midrule
& $N_h$ & $L^1$ Error & Order & $L^\infty$ Error & Order & CPU(s) \\
\midrule
& 160 & $1.27 \times 10^{-5}$ & - & $4.91 \times 10^{-5}$ & - & 1.83 \\
& 320 & $1.59 \times 10^{-6}$ & 3.00 & $6.14 \times 10^{-6}$ & 3.00 & 11.74\\
& 640 & $1.98 \times 10^{-7}$ & 3.00 & $7.68 \times 10^{-7}$ & 3.00 & 80.38\\
& 1280 & $2.47 \times 10^{-8}$ & 3.00 & $9.60 \times 10^{-8}$ & 3.00 & 748.34\\
\midrule
\multicolumn{7}{c}{Sparse grid} \\
\midrule
$N_r$ & $N_h$ & $L^1$ Error & Order & $L^\infty$ Error & Order & CPU(s) \\
\midrule
20 & 160 & $4.56 \times 10^{-5}$ & - & $5.21 \times 10^{-4}$ & - & 1.51 \\
40 & 320 & $2.11 \times 10^{-6}$ & 4.43 & $1.63 \times 10^{-4}$ & 1.67 & 6.08 \\
80 & 640  & $2.73 \times 10^{-7}$ & 2.95 & $1.35 \times 10^{-5}$ & 3.60 & 22.63 \\
160 & 1280 & $3.02 \times 10^{-8}$ & 3.18 & $1.56 \times 10^{-6}$ & 3.12 & 159.19 \\
\bottomrule[1.5pt]
\end{tabular}
\caption{Example 1, a linear problem with smooth solution. RK FPFS scheme with the third order linear upwind
approximations, comparison of numerical errors and CPU times for computations on single-grid and sparse-grid.
Third order Lagrange interpolation for prolongation is employed in sparse-grid computations.
$N_r$: number of cells in each spatial direction of a root grid.
CPU: CPU time for a complete simulation. CPU time unit: seconds.}
\label{tab:smooththird}
\end{table}


\paragraph{Example 2 (A nonlinear problem with smooth solution).}
We solve the 2D Eikonal equation \eqref{eq:eikonal}  with the right hand side function
$$
f(x,y) = \frac{\pi}{2} \sqrt{\sin^2 \left(\pi + \frac{\pi}{2}x \right) +\sin^2 \left(\pi + \frac{\pi}{2}y \right)  },
$$
and the source point $\Gamma = (0,0)$. The computational domain $\Omega  = [-1,1]^2$. The exact solution of the problem is
$$
\phi(x,y) = \cos \left( \pi + \frac{\pi}{2} x\right) + \cos \left( \pi + \frac{\pi}{2} y\right).
$$
We use this example to verify that the proposed sparse grid RK FPFS-WENO scheme can achieve the desired
accuracy order for a nonlinear problem with smooth solution. The sparse grid RK FPFS-WENO scheme with $\gamma = 0.4$ and the third order WENO approximations to the derivatives is applied.
Both the third order Lagrange interpolation and the third order WENO interpolation are used for prolongation in sparse-grid computations. We perform simulations on both sparse grids and the corresponding single grids, and compare their results. The $L^1$ errors, $L^\infty$ errors and their numerical accuracy orders and CPU times are reported in Table~\ref{tab:source}. It is observed that third order accuracy is obtained for all cases, including sparse-grid computations with Lagrange or WENO prolongation and the corresponding single-grid ones, along with the mesh refinement. Comparing the numerical errors of sparse-grid computations and the corresponding single-grid ones, similar as Example 1 we observe that
their $L^1$ errors are comparable while sparse-grid computations with WENO prolongation have slight larger errors. The $L^\infty$ errors of sparse-grid computations are larger than the corresponding single-grid computations. In terms of computational efficiency, on refined meshes we see around $65\% \sim 70\%$ CPU time saved for simulations on sparse grids vs single grids, for this nonlinear example. We also notice that on relatively coarse mesh (e.g. $N_h=160$ here), it takes more CPU time for sparse-grid computation than the corresponding single-grid one, due to quite different iteration history on different semi-coarsened sparse grids.

\begin{table}
\centering
\begin{tabular}{ccccccc}
\toprule[1.5pt]
\multicolumn{7}{c}{Single grid} \\
\midrule
& $N_h$ & $L^1$ Error & Order & $L^\infty$ Error & Order & CPU(s)\\
\midrule
& 160 & $1.05 \times 10^{-6}$ & - & $1.78 \times 10^{-6}$ & - & 5.55\\
& 320 & $1.11 \times 10^{-7}$ & 3.24 & $1.71 \times 10^{-7}$ & 3.38 & 30.59\\
&  640 & $1.37 \times 10^{-8}$ & 3.02 & $2.10 \times 10^{-8}$ & 3.30 & 297.53\\
& 1280 & $1.71 \times 10^{-9}$ & 3.00 & $2.61 \times 10^{-9}$ & 3.00 & 1,401.50\\
\midrule
\multicolumn{7}{c}{Sparse grid, Lagrange interpolation} \\
\midrule
$N_r$ & $N_h$ & $L^1$ Error & Order & $L^\infty$ Error & Order & CPU(s)\\
\midrule
20 & 160 & $3.28 \times 10^{-6}$ & - & $1.74 \times 10^{-5}$ & - & 11.50 \\
40 & 320 & $2.70 \times 10^{-7}$ & 3.60 & $2.96 \times 10^{-6}$ & 2.55 & 28.28\\
80 & 640 & $2.34 \times 10^{-8}$ & 3.53 & $4.20 \times 10^{-7}$ & 2.82 & 94.45 \\
160 & 1280 & $2.27 \times 10^{-9}$ & 3.36 & $5.55 \times 10^{-8}$ & 2.92 & 485.06\\
\midrule
\multicolumn{7}{c}{Sparse grid, WENO interpolation} \\
\midrule
$N_r$ & $N_h$ & $L^1$ Error & Order & $L^\infty$ Error & Order & CPU(s) \\
\midrule
20 & 160  & $8.60 \times 10^{-6}$ & - & $4.72 \times 10^{-3}$ & - & 11.74\\
40 & 320 & $7.65 \times 10^{-7}$ & 3.49 & $1.21 \times 10^{-3}$ & 1.96 & 29.57\\
80 & 640 & $6.00 \times 10^{-8}$ & 3.67 & $2.72 \times 10^{-4}$ & 2.16 & 87.74\\
160 & 1280 & $4.17 \times 10^{-9}$ & 3.84 & $2.00 \times 10^{-5}$ & 3.76 & 494.72\\
\bottomrule[1.5pt]
\end{tabular}
\caption{Example 2, a nonlinear problem with smooth solution. RK FPFS-WENO scheme, comparison of numerical errors and CPU times for computations on single-grid and sparse-grid.
Both third order Lagrange interpolation and WENO interpolation for prolongation are employed in sparse-grid computations. $N_r$: number of cells in each spatial direction of a root grid.
CPU: CPU time for a complete simulation. CPU time unit: seconds.
}
\label{tab:source}
\end{table}


\bigskip
\bigskip
\noindent In the following, we apply the sparse grid scheme to examples with non-smooth solutions to show its nonlinear stability and computational efficiency.

\begin{figure}
	\centering
	\begin{subfigure}[b]{0.45\textwidth}
		\includegraphics[width=\textwidth]{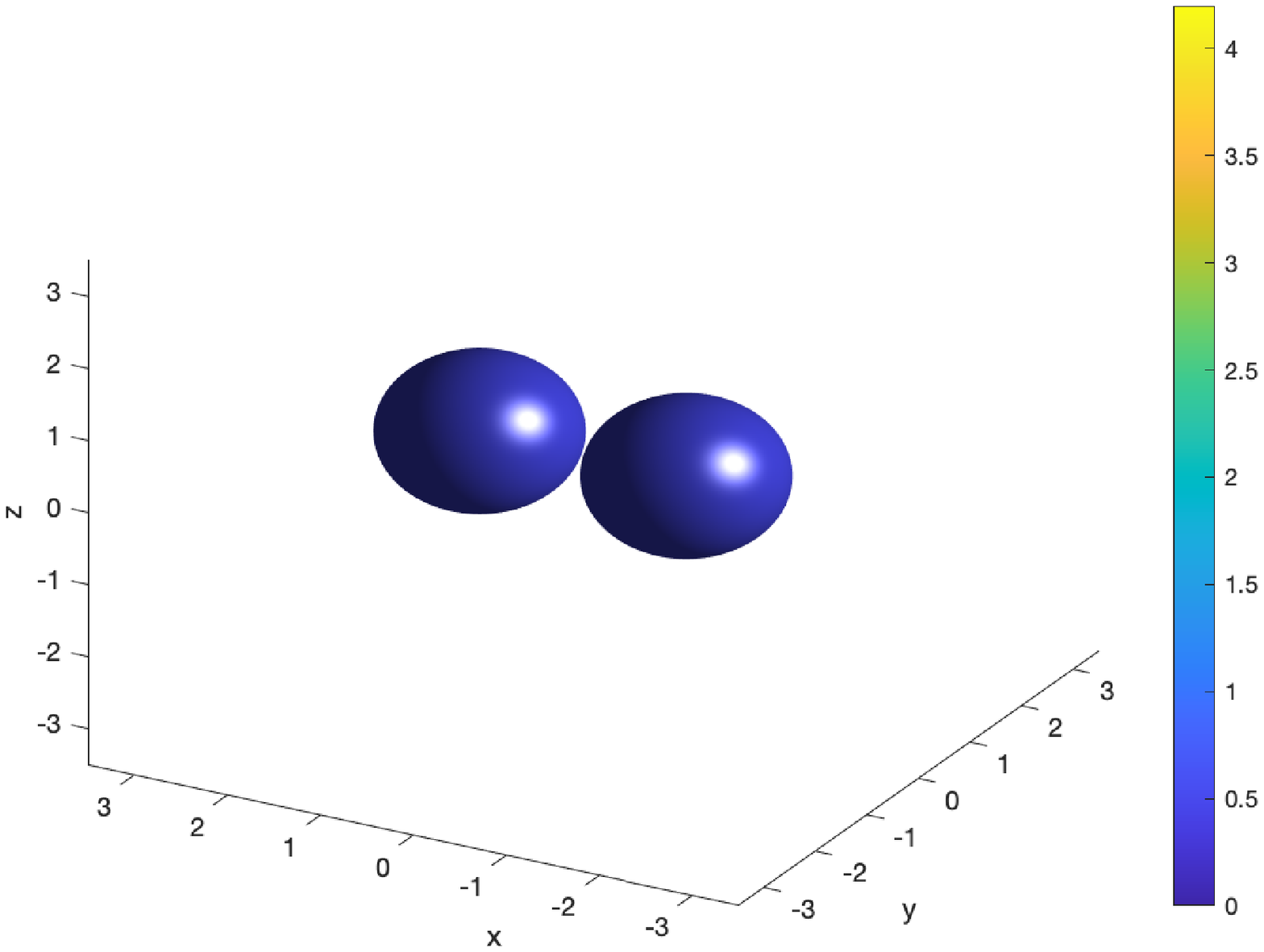}
		\caption{$\phi = 0.5$}
	\end{subfigure}
	\begin{subfigure}[b]{0.45\textwidth}
		\includegraphics[width=\textwidth]{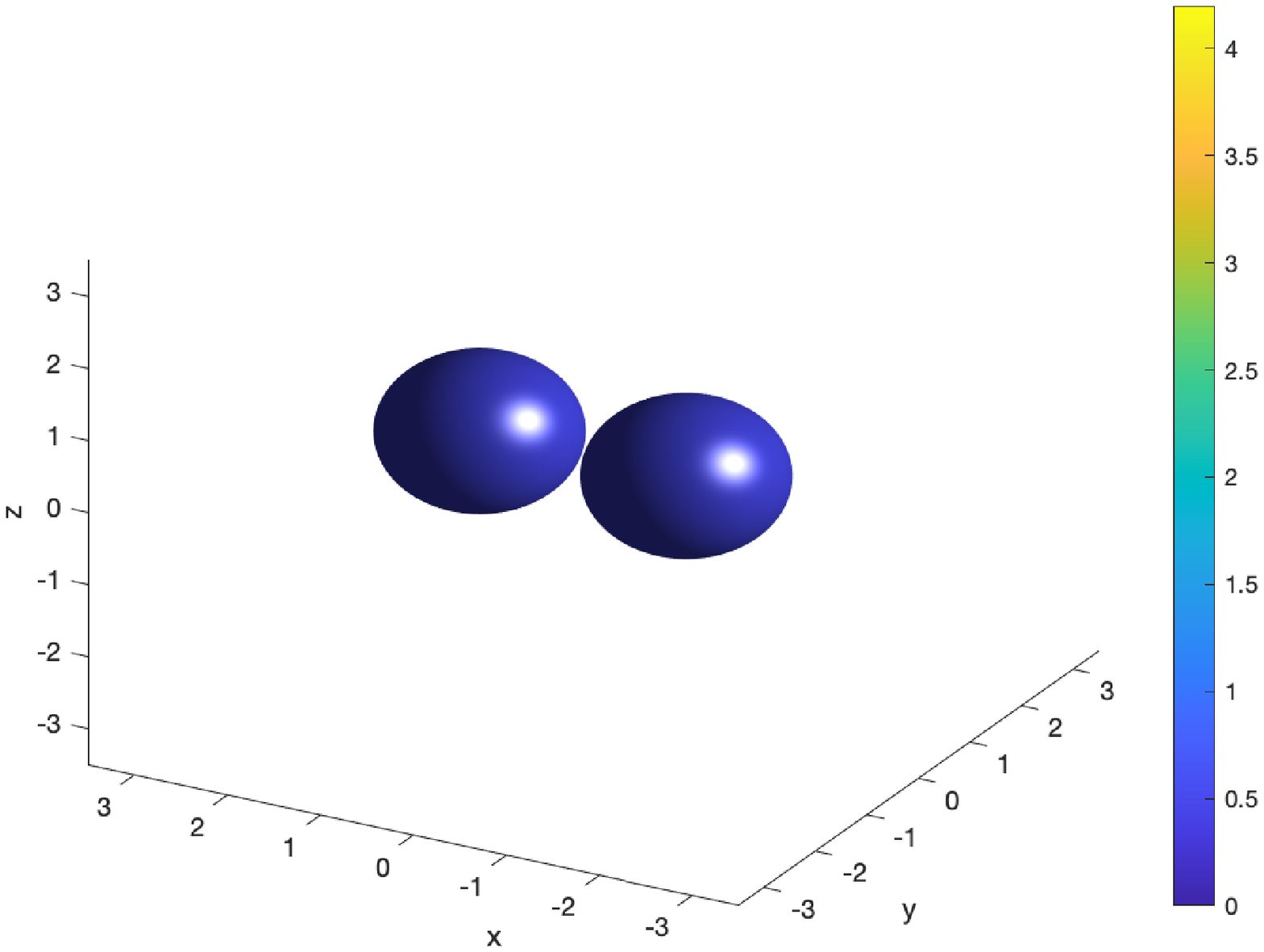}
		\caption{$\phi = 0.5$}
	\end{subfigure}

	\begin{subfigure}[b]{0.45\textwidth}
		\includegraphics[width=\textwidth]{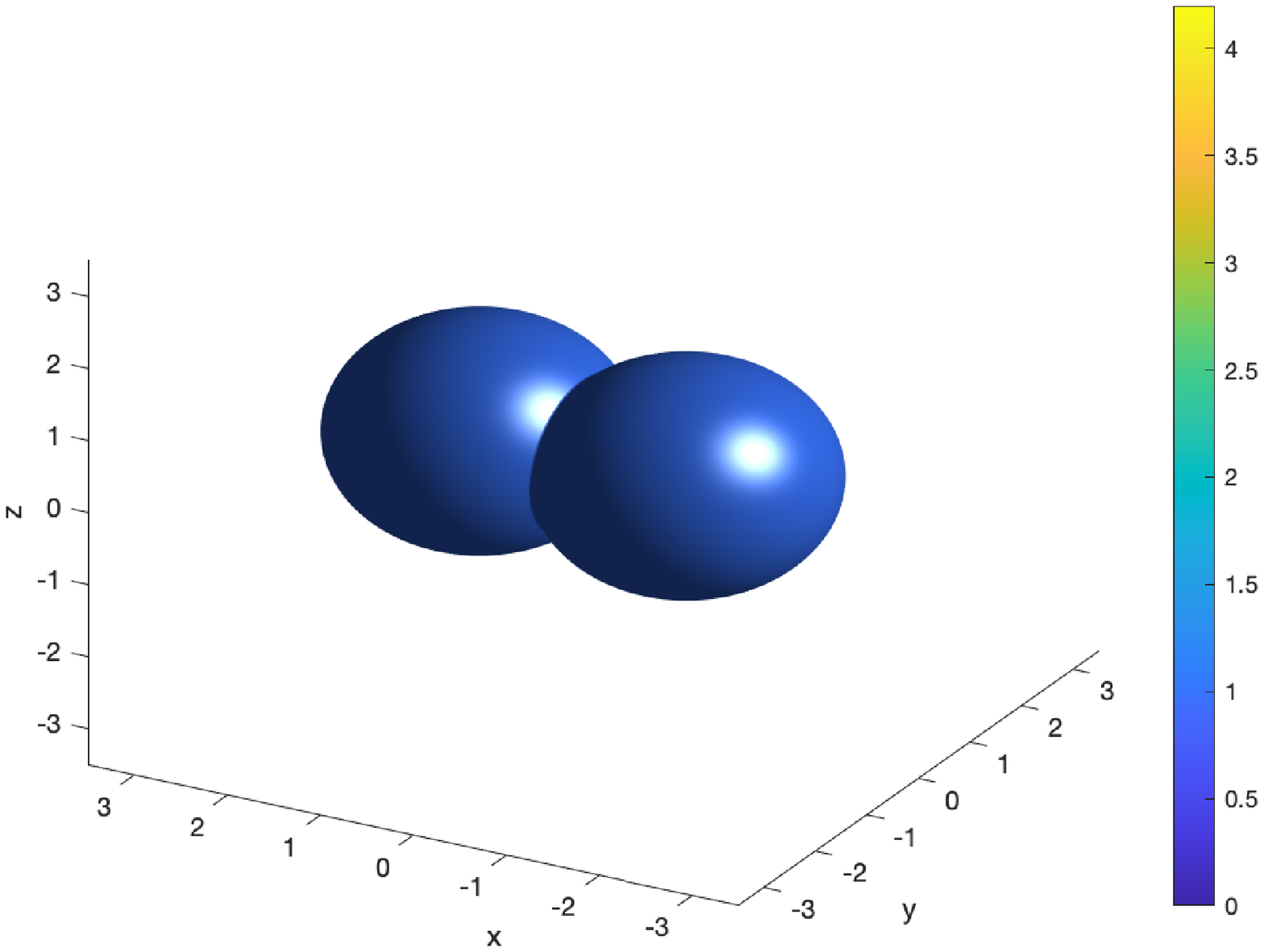}
		\caption{$\phi = 1$}
	\end{subfigure}
	\begin{subfigure}[b]{0.45\textwidth}
		\includegraphics[width=\textwidth]{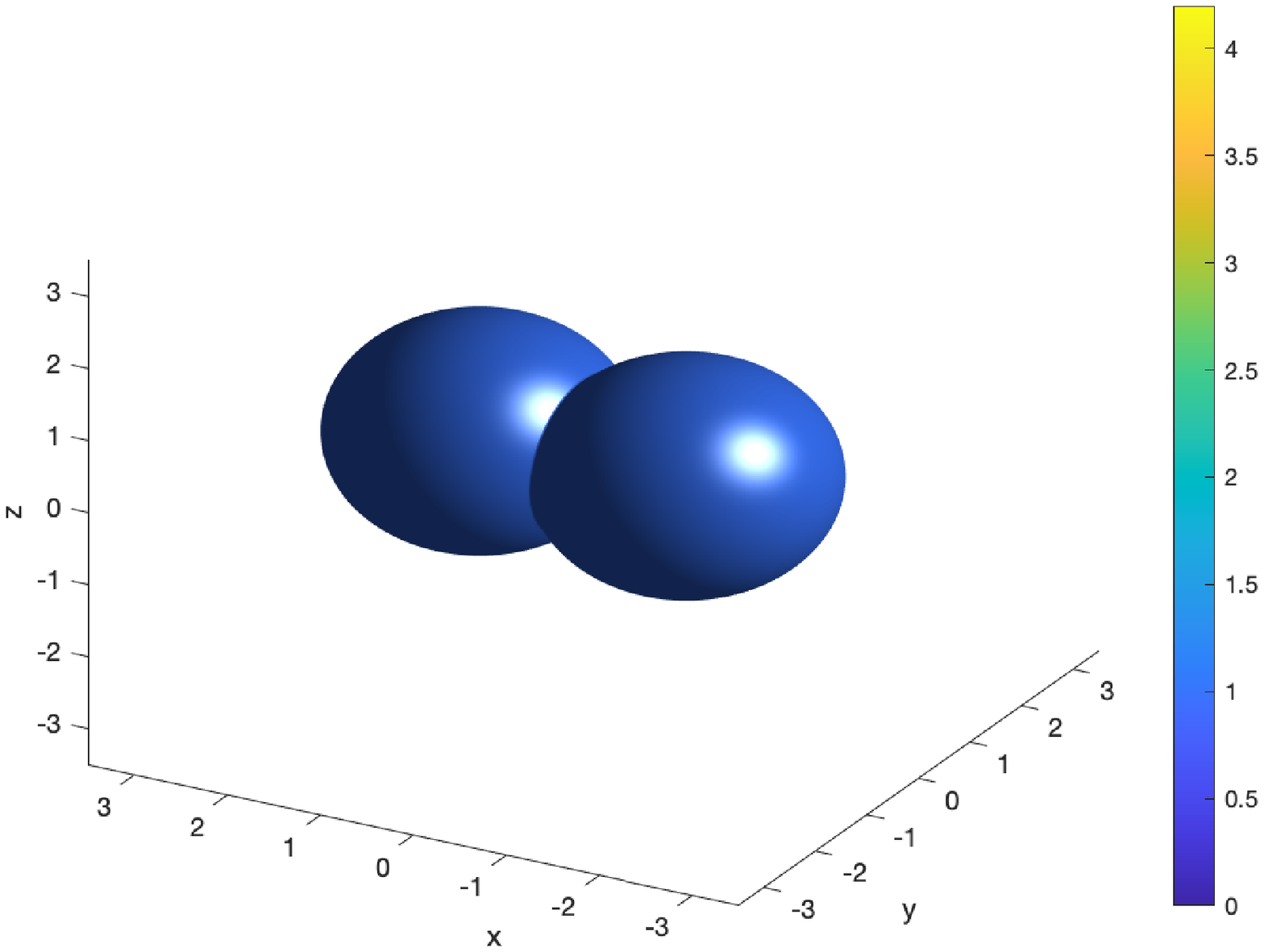}
		\caption{$\phi = 1$}
	\end{subfigure}

	\begin{subfigure}[b]{0.45\textwidth}
		\includegraphics[width=\textwidth]{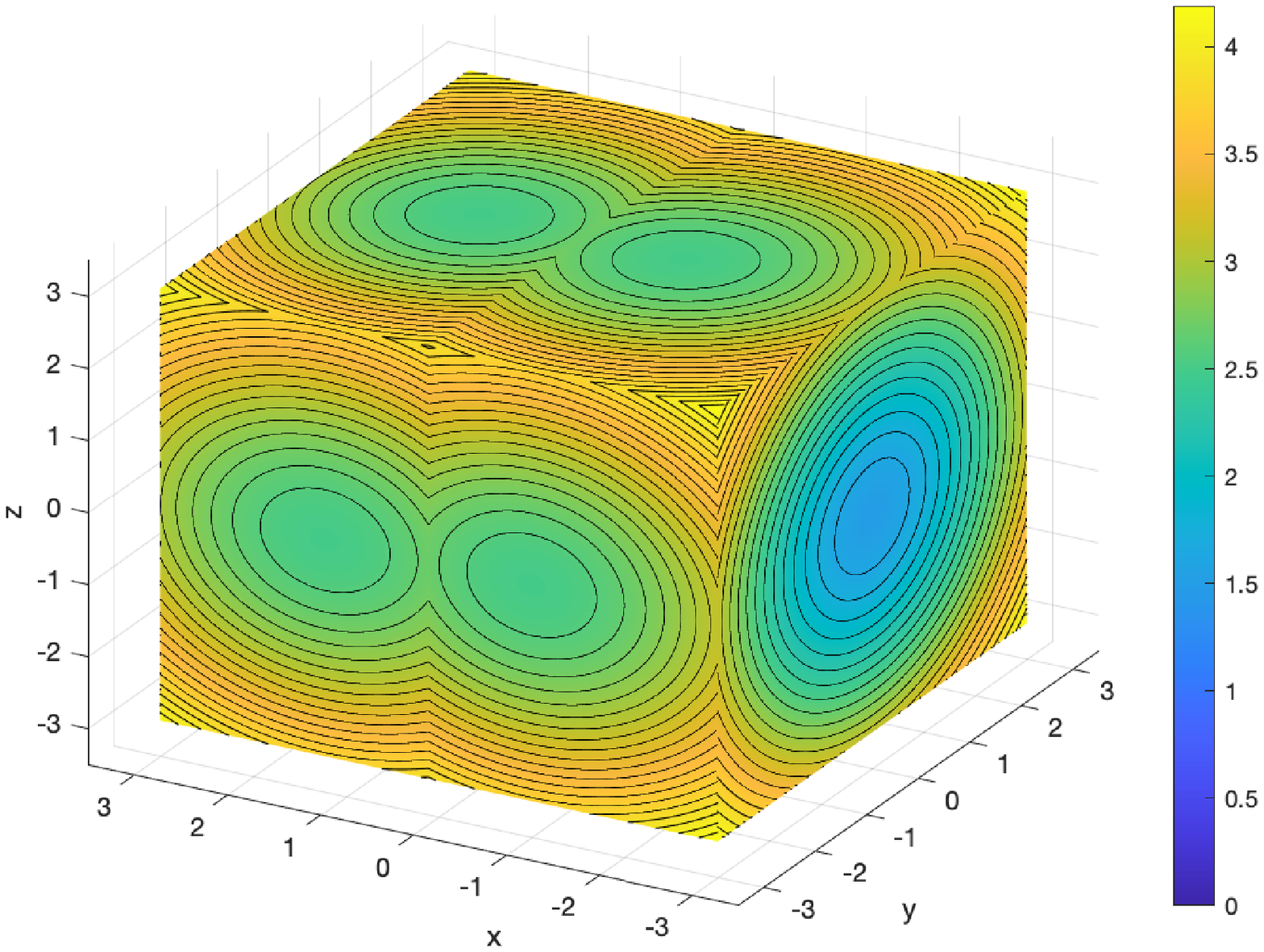}
		\caption{Surface contours}
	\end{subfigure}
	\begin{subfigure}[b]{0.45\textwidth}
		\includegraphics[width=\textwidth]{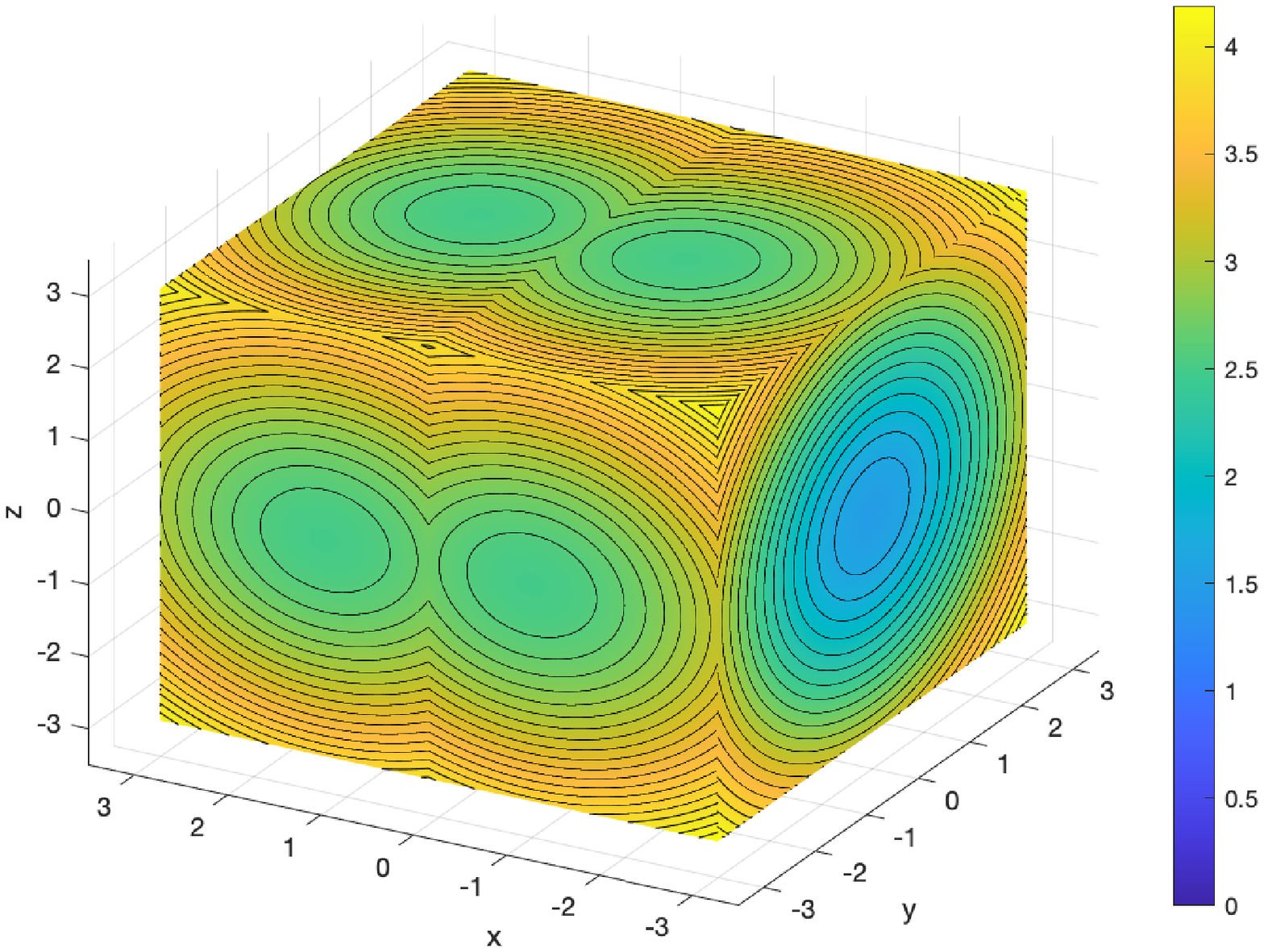}
		\caption{Surface contours}
	\end{subfigure}
	\caption{\footnotesize{Example 3, numerical solutions of the two-sphere problem by the RK FPFS-WENO scheme on sparse grids ($N_r = 80$ for root grid, finest level $N_L=3$ in the sparse-grid computation) and the corresponding $640\times640\times640$ single grid, using the third order WENO interpolation for prolongation in the sparse-grid combination. (a), (c), (e): single-grid result; (b), (d), (f): sparse-grid result; (a), (b): the contour plots for $\phi=0.5$; (c), (d): the contour plots for $\phi=1$; (e), (f): the contour plots for the whole surface.}}
	\label{fig:2sphere}
\end{figure}

\paragraph{Example 3 (Two-sphere problem).}
We solve the 3D Eikonal equation \eqref{eq:eikonal}  with $f(x,y,z) = 1$ on the computational domain $\Omega = [-3,3] ^3$. $\Gamma$ are two spheres of equal radius 0.5 centered at $(-1,0,0)$ and $(\sqrt{1.5},0,0)$. The exact solution of the problem is the distance function to $\Gamma$: $\phi(x,y,z) = \min(d_1,d_2)$, where
\begin{align*}
d_1 &= \left| \sqrt{(x+1)^2 + y^2 + z^2} - 0.5 \right|, \\
d_2 &= \left| \sqrt{(x-\sqrt{1.5})^2 + y^2 + z^2} - 0.5 \right|.
\end{align*}
The solution of the problem is non-smooth. Singularities exist in the centers of each sphere and the plane that is equidistant from both spheres. The sparse grid RK FPFS-WENO scheme with $\gamma = 0.8$ and the third order WENO approximations to the derivatives is applied. The third order WENO interpolation is used for prolongation in sparse-grid computations. We perform simulations on both sparse grids with $N_r = 80, N_L=3$ and the corresponding $640\times 640\times 640$ single grid, and compare their results.
The results are shown in Fig.~\ref{fig:2sphere}. We observe that the numerical solutions by the sparse grid RK FPFS-WENO scheme and its corresponding single-grid simulation are comparable. The nonlinear stability and high resolution properties of the RK FPFS-WENO scheme for resolving the non-smooth solution are preserved well in the sparse-grid simulation.
We record the simulation CPU time costs to compare the computational efficiency. It takes $42,674.31$ seconds of CPU time to complete the simulation in the sparse-grid computation, while $458,311.48$ seconds of CPU time are needed for finishing the simulation in the corresponding single-grid computation. About $91\%$ CPU time is saved by performing the RK FPFS-WENO simulation on the sparse grids here.

\begin{figure}
	\centering
	\begin{subfigure}[b]{0.45\textwidth}
		\includegraphics[width=\textwidth]{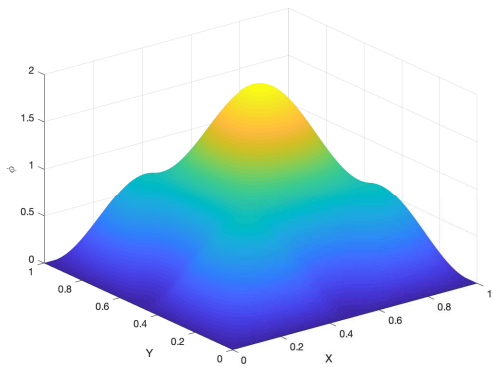}
		\caption{single-grid result}
	\end{subfigure}
	\begin{subfigure}[b]{0.45\textwidth}
		\includegraphics[width=\textwidth]{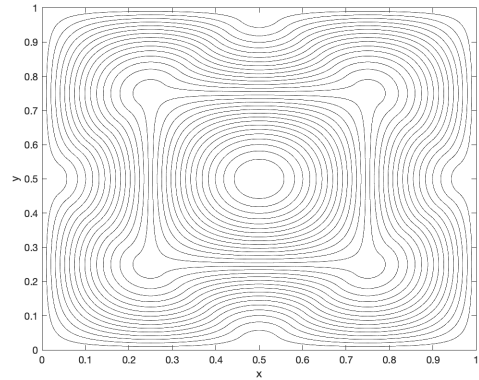}
		\caption{single-grid result}
	\end{subfigure}

	\begin{subfigure}[b]{0.45\textwidth}
		\includegraphics[width=\textwidth]{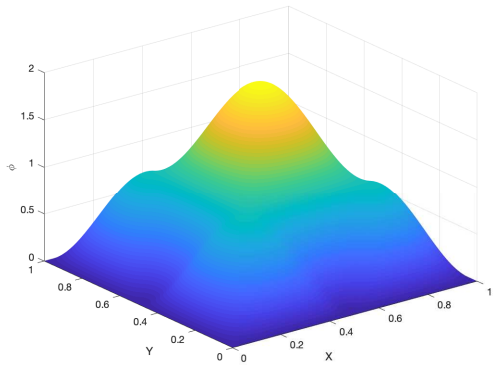}
		\caption{sparse-grid result}
	\end{subfigure}
	\begin{subfigure}[b]{0.45\textwidth}
		\includegraphics[width=\textwidth]{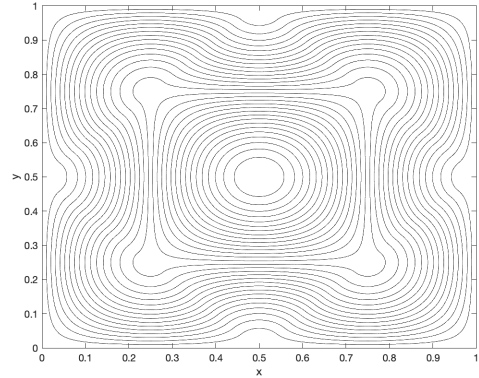}
		\caption{sparse-grid result}
	\end{subfigure}

	\caption{\footnotesize{Example 4, numerical solutions of the shape-from-shading problem by the RK FPFS-WENO scheme on sparse grids ($N_r = 160$ for root grid, finest level $N_L=3$ in the sparse-grid computation) and the corresponding $1280\times1280$ single grid, using the third order WENO interpolation for prolongation in the sparse-grid combination. (a), (b): single-grid result; (c), (d): sparse-grid result; (a), (c): three-dimensional view of the solutions; (b), (d): the contour plots, 30 equally spaced contour lines from $\phi = 0$ to $\phi = 2$.}}
	\label{fig:shape2}
\end{figure}

\paragraph{Example 4 (Shape-from-shading).}
We solve the Eikonal equation \eqref{eq:eikonal}  with the right hand side function
$$
f(x,y) = 2\pi  \sqrt{[\cos(2\pi x)\sin(2\pi y)]^2 + [\sin(2\pi x)\cos(2\pi y)]^2}.
$$
The computational domain $\Omega = [0,1] \times [0,1]$. $\phi(x,y) = 0$  is prescribed at the boundary $\partial\Omega$ of the unit square. The boundary region $\Gamma = \{ (\frac{1}{4},\frac{1}{4}), (\frac{3}{4},\frac{3}{4}), (\frac{1}{4},\frac{3}{4}), (\frac{3}{4},\frac{1}{4}), (\frac{1}{2},\frac{1}{2}) \}\cup\partial\Omega$, consisting of five isolated points and $\partial\Omega$. The values at these five isolated points are specified as
$$
g\left(\frac{1}{4},\frac{1}{4} \right) = g\left(\frac{3}{4},\frac{3}{4} \right) =  g\left(\frac{1}{4},\frac{3}{4} \right) = g\left(\frac{3}{4},\frac{1}{4} \right) = 1, \qquad g\left(\frac{1}{2},\frac{1}{2} \right)  = 2.
$$
The exact solution of the problem is
$$
\phi(x,y) = \left\{ \begin{matrix}
\max(|\sin(2\pi x) \sin(2\pi y)|, 1+ \cos(2\pi x)\cos(2\pi  y)), \\
\quad \text{if } |x+y-1| < \frac{1}{2} \text{ and } |x-y| < \frac{1}{2}; \\
|\sin(2\pi x) \sin(2\pi y)|, \qquad \text{otherwise},
\end{matrix} \right.
$$
which is {\em not} smooth.
Actually the solution of this problem is the shape function, which has the brightness $I(x,y) = 1/\sqrt{1+f(x,y)^2}$ under vertical lighting. Details about this problem can be found in \cite{RTNum}.
The sparse grid RK FPFS-WENO scheme with
$\gamma = 0.4$ and the third order WENO approximations to the derivatives is applied. The third order WENO interpolation is used for prolongation in sparse-grid computations. Simulations are carried out on both sparse grids with $N_r = 160, N_L=3$ and the corresponding $1280\times 1280$ single grid, to compare their results.
The results are reported in Fig.~\ref{fig:shape2}. It is observed that the numerical solutions by the sparse grid RK FPFS-WENO scheme and its corresponding single-grid computation are comparable. As the previous example, the nonlinear stability and high resolution properties of the RK FPFS-WENO scheme for resolving the non-smooth solution of this example are preserved well in the sparse-grid computation.
Again, we record the simulation CPU time costs to compare their computational efficiency. It takes $670.86$  seconds of CPU time to complete the simulation in the sparse-grid computation, while $1,510.27$ seconds of CPU time are needed for finishing the simulation in the corresponding single-grid computation. About $56\%$ CPU time is saved by carrying out the RK FPFS-WENO simulation on the sparse grids in this example.

\paragraph{Example 5 (Voronoi diagram problem).}

We consider a Voronoi diagram problem as in \cite{Auren1991,OBSC2000}. Given a set of points (called generators) in a domain, the Voronoi diagram divides the domain into regions in which all points inside the region are closest to the generator of that region than any other  generators. This kind of problems have applications in many fields, including engineering, natural sciences, geometry, humanities, etc, for example, dividing a map into response regions for local fire stations. An essential part for solving a Voronoi diagram problem is to compute the minimum travel time to the closest generator by solving the Eikonal equation \eqref{eq:eikonal}. Here we solve both a 2D case and a 3D case.

\begin{figure}
	\centering
	\begin{subfigure}[b]{0.45\textwidth}
		\includegraphics[width=\textwidth]{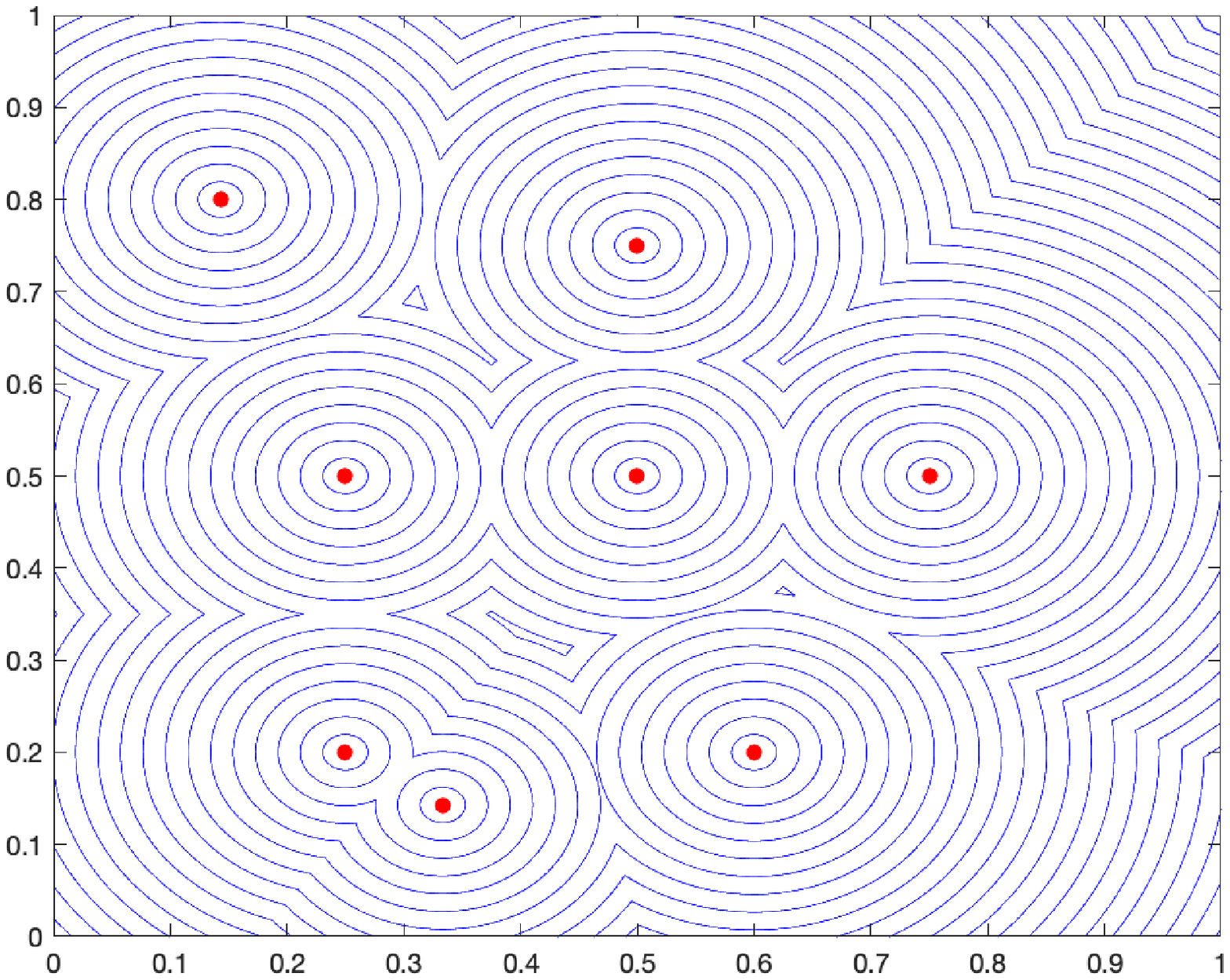}
		\caption{single-grid result}
	\end{subfigure}
	\begin{subfigure}[b]{0.45\textwidth}
		\includegraphics[width=\textwidth]{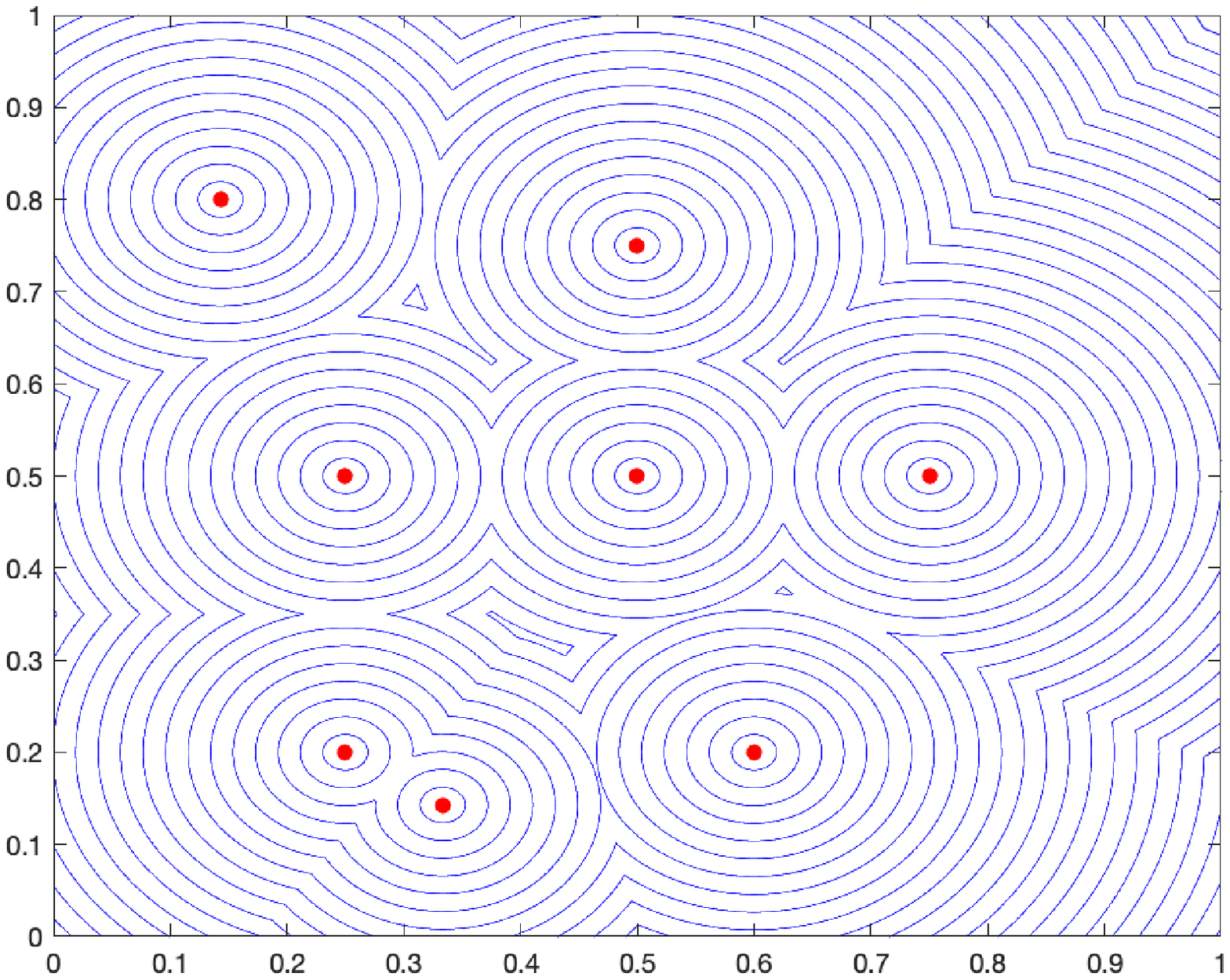}
		\caption{sparse-grid result}
	\end{subfigure}
	\caption{\footnotesize{Example 5, Case 1, numerical solutions of the 2D Voronoi diagram problem by the RK FPFS-WENO scheme on sparse grids ($N_r = 160$ for root grid, finest level $N_L=3$ in the sparse-grid computation) and the corresponding $1280\times1280$ single grid, using the third order WENO interpolation for prolongation in the sparse-grid combination. The contour plots,  30 equally spaced contour lines from $\phi = 0$ to $\phi = 0.5589$. Red points are the generators. (a): single-grid result; (b): sparse-grid result.}}
	\label{fig:2dvor}
\end{figure}

\textbf{Case 1 (2D).} We solve the Eikonal equation \eqref{eq:eikonal} with $f(x,y) = 1$. The computational domain  $\Omega = [0,1]^2$. $\phi(x,y) = 0$ is prescribed at the points (the generators):
$$
\Gamma = \left\{ \left( \frac{1}{4}, \frac{1}{5} \right), \left( \frac{1}{3}, \frac{1}{7} \right), \left( \frac{3}{5}, \frac{1}{5} \right), \left( \frac{3}{4}, \frac{1}{2} \right),  \left( \frac{1}{2}, \frac{3}{4}  \right), \left( \frac{1}{4}, \frac{1}{2} \right), \left( \frac{1}{7}, \frac{4}{5} \right), \left( \frac{1}{2}, \frac{1}{2} \right) \right\}.
$$
The exact solution of the problem is the distance function to $\Gamma$, and it is {\em not} smooth.
The sparse grid RK FPFS-WENO scheme with $\gamma = 0.8$ and the third order WENO approximations to the derivatives is applied. The third order WENO interpolation is used for prolongation in sparse-grid computations. Simulations are performed on both sparse grids with $N_r = 160, N_L=3$ and the corresponding $1280\times 1280$ single grid, to compare their results.
The results are reported in Fig.~\ref{fig:2dvor}. Again, we observe that the numerical solutions by the sparse grid RK FPFS-WENO scheme and its corresponding single-grid computation are comparable, and the nonlinear stability and high resolution properties of the RK FPFS-WENO scheme for resolving the non-smooth solution in this example are preserved well in the sparse-grid simulation. About computational efficiency, it takes $634.42$  seconds of CPU time to complete the simulation in the sparse-grid computation, while $1,591.41$ seconds of CPU time are needed for finishing the simulation in the corresponding single-grid computation. About $60\%$ CPU time is saved by performing the RK FPFS-WENO simulation on the sparse grids in this problem.

\begin{figure}
	\centering
	\begin{subfigure}[b]{0.45\textwidth}
		\includegraphics[width=\textwidth]{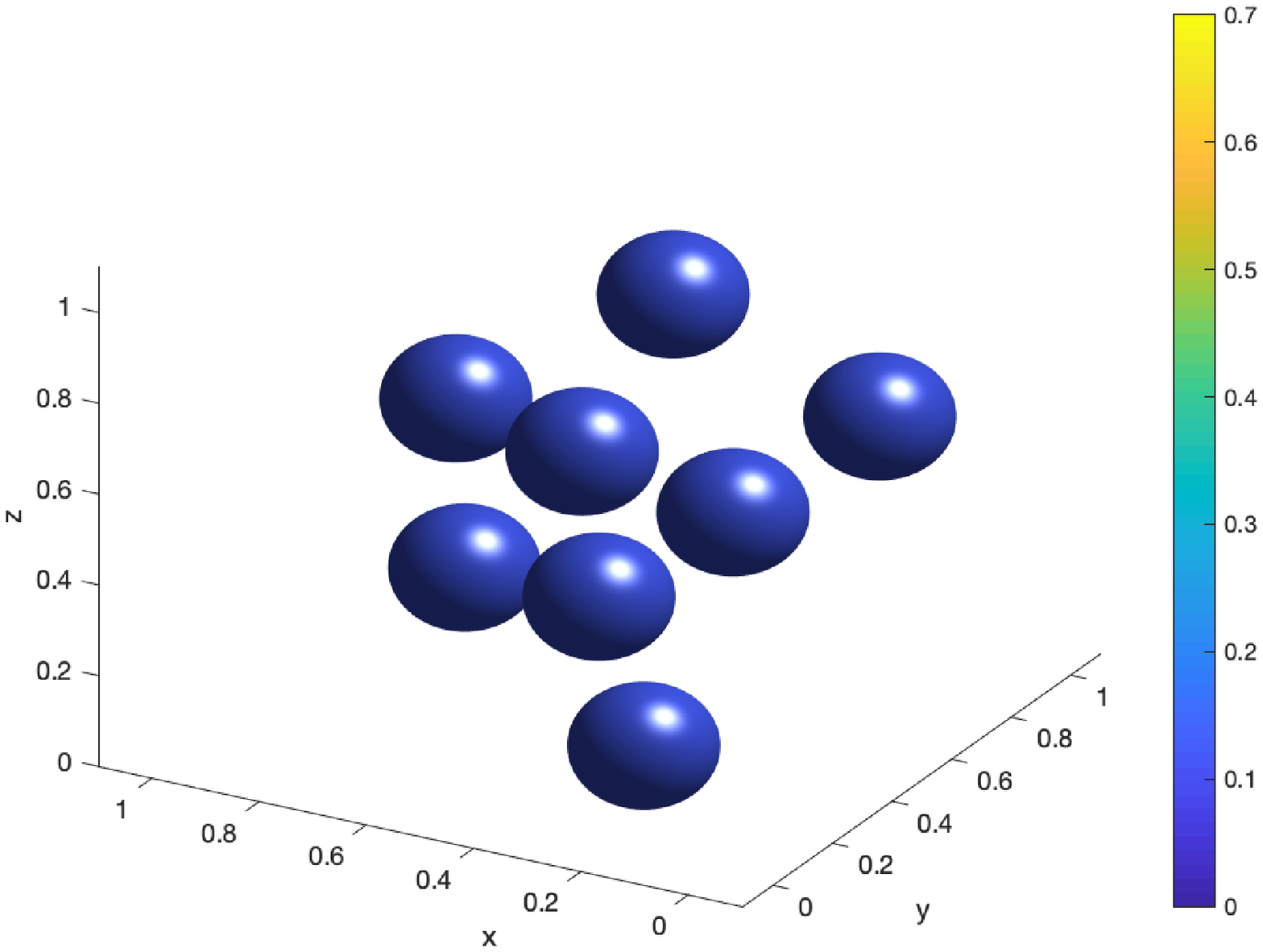}
		\caption{$\phi = 0.125$}
	\end{subfigure}
	\begin{subfigure}[b]{0.45\textwidth}
		\includegraphics[width=\textwidth]{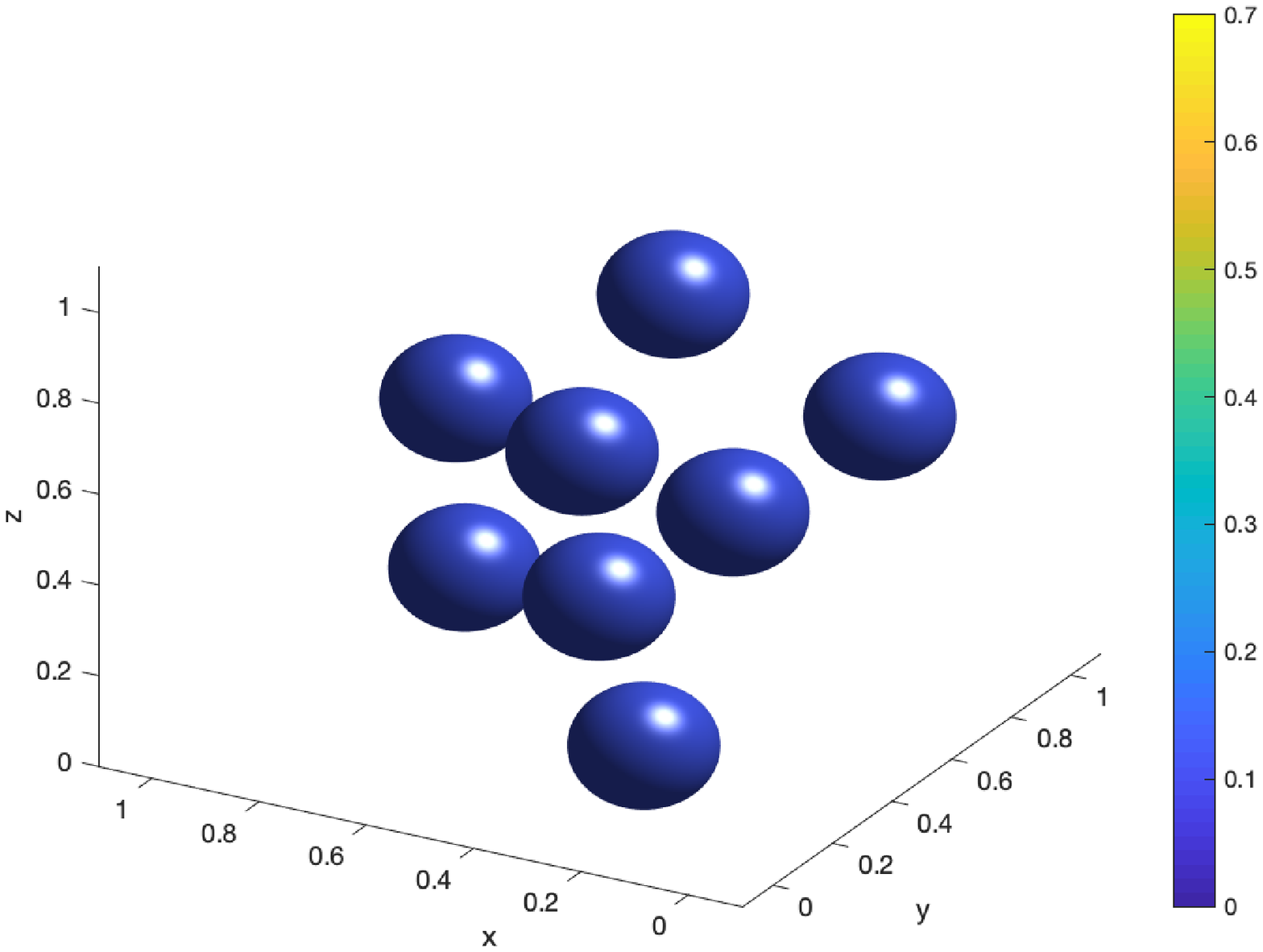}
		\caption{$\phi = 0.125$}
	\end{subfigure}

	\begin{subfigure}[b]{0.45\textwidth}
		\includegraphics[width=\textwidth]{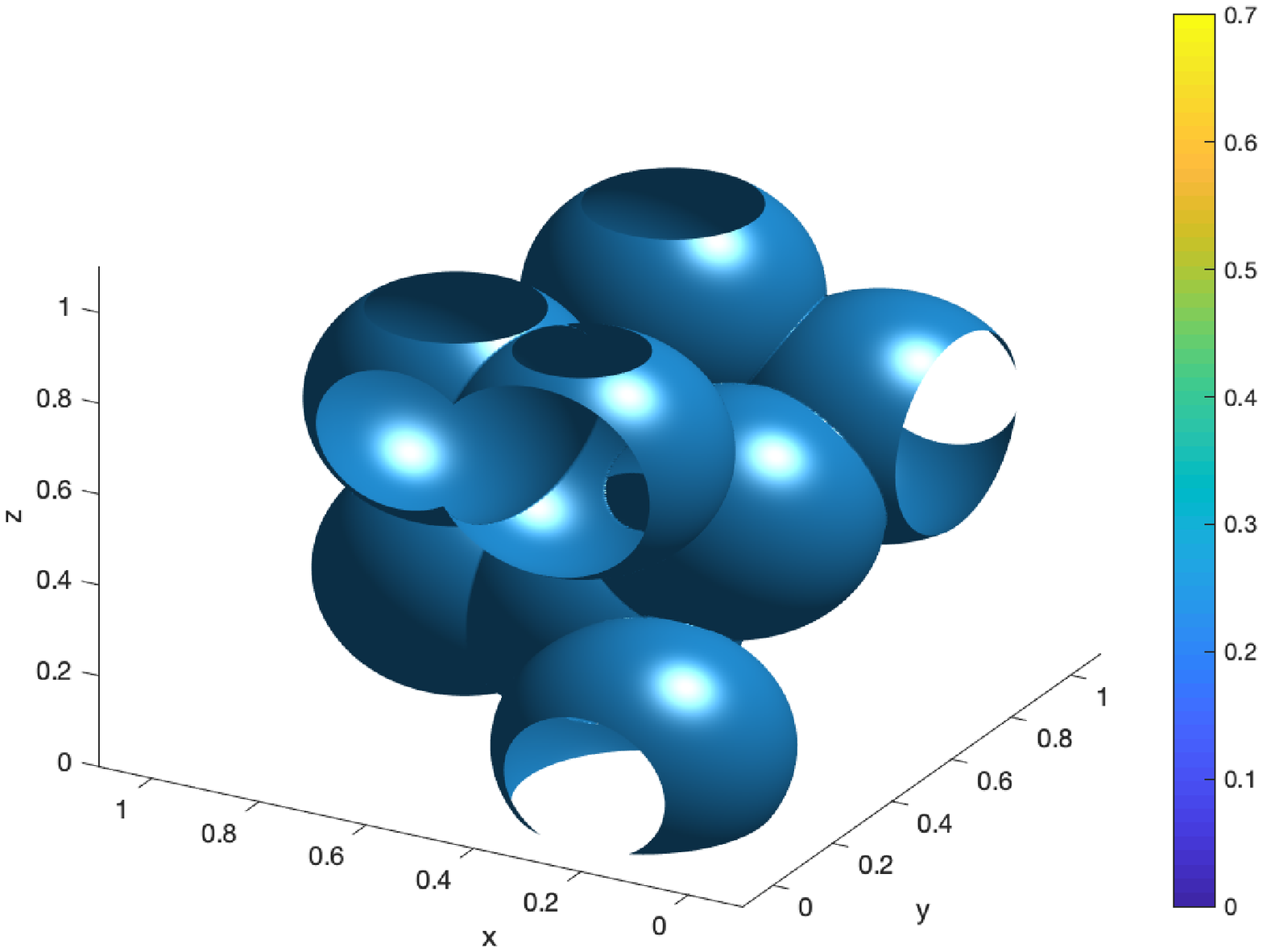}
		\caption{$\phi = 0.25$}
	\end{subfigure}
	\begin{subfigure}[b]{0.45\textwidth}
		\includegraphics[width=\textwidth]{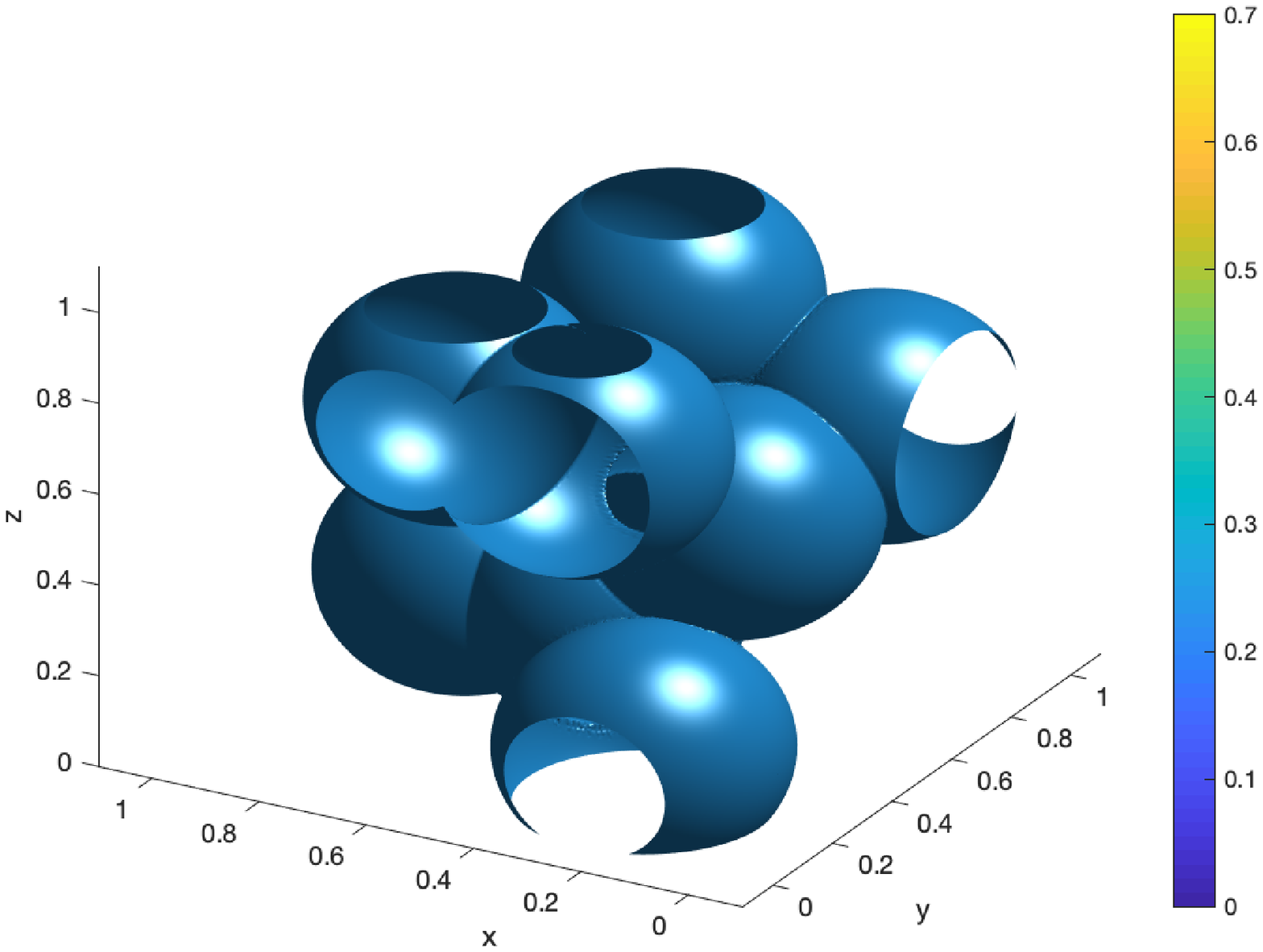}
		\caption{$\phi = 0.25$}
	\end{subfigure}

	\begin{subfigure}[b]{0.45\textwidth}
		\includegraphics[width=\textwidth]{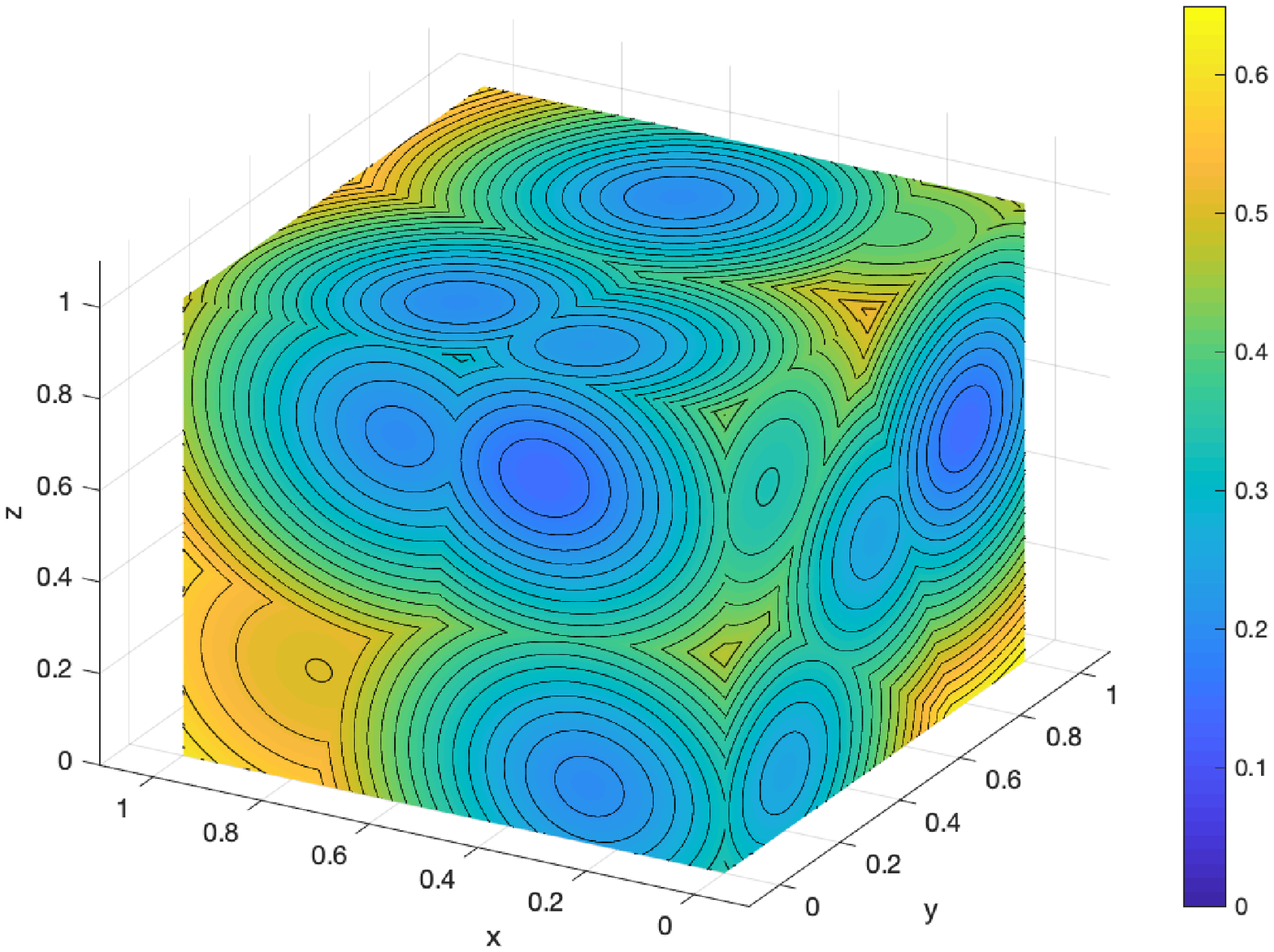}
		\caption{Surface contour}
	\end{subfigure}
	\begin{subfigure}[b]{0.45\textwidth}
		\includegraphics[width=\textwidth]{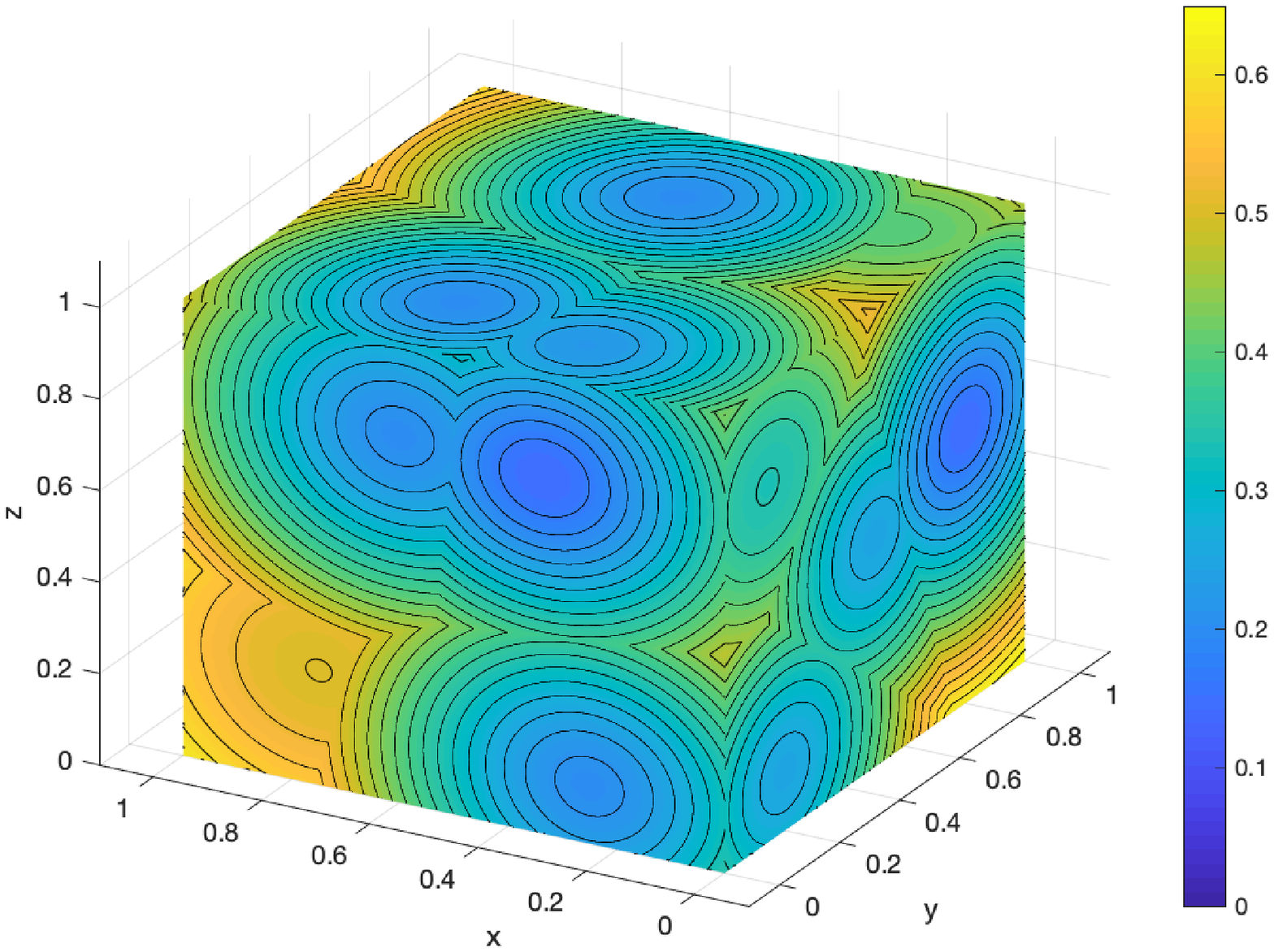}
		\caption{Surface contour}
	\end{subfigure}
	\caption{\footnotesize{Example 5, Case 2, numerical solutions of the 3D Voronoi diagram problem by the RK FPFS-WENO scheme on sparse grids ($N_r = 80$ for root grid, finest level $N_L=3$ in the sparse-grid computation) and the corresponding $640\times640\times640$ single grid, using the third order WENO interpolation for prolongation in the sparse-grid combination. (a), (c), (e): single-grid result; (b), (d), (f): sparse-grid result; (a), (b): the contour plots for $\phi=0.125$; (c), (d): the contour plots for $\phi=0.25$; (e), (f): the contour plots for the whole surface.}}
	\label{fig:3dvor}
\end{figure}

\textbf{Case 2 (3D).} Now we solve the 3D case, the Eikonal equation \eqref{eq:eikonal} with $f(x,y,z) = 1$. The computational domain  $\Omega = [0,1]^3$. $\phi(x,y,z) = 0$ is specified at the following generators
\begin{align*}
\Gamma &= \left\{ \left( \frac{1}{4}, \frac{1}{5}, \frac{1}{8}\right), \left( \frac{1}{3}, \frac{1}{7}, \frac{7}{9}\right), \left( \frac{3}{5}, \frac{1}{5}, \frac{4}{5}\right), \left( \frac{3}{4}, \frac{1}{2}, \frac{1}{4}\right), \right. \\
&\qquad \left. \left( \frac{1}{2}, \frac{3}{4}, \frac{4}{5}\right), \left( \frac{1}{4}, \frac{1}{2}, \frac{1}{2}\right), \left( \frac{1}{7}, \frac{4}{5}, \frac{3}{5}\right), \left( \frac{1}{2}, \frac{1}{2}, \frac{1}{4}\right) \right\}.
\end{align*}
The exact solution of the problem is the distance function to $\Gamma$ in this 3D domain, and it is {\em not} a smooth function.
The sparse grid RK FPFS-WENO scheme with $\gamma = 0.8$ and the third order WENO approximations to the derivatives is applied in solving this problem. The third order WENO interpolation is used for prolongation in sparse-grid computations. Simulations are performed on both sparse grids with $N_r = 80, N_L=3$ and the corresponding $640\times 640\times 640$ single grid, to compare their results.
The simulation results are presented in Fig.~\ref{fig:3dvor}, which shows that the numerical solutions by the sparse grid RK FPFS-WENO scheme and its corresponding single-grid computation are comparable, and the nonlinear stability and high resolution properties of the RK FPFS-WENO scheme for resolving the non-smooth solution in this 3D example are preserved well in the sparse-grid simulation. In terms of computational efficiency, it takes $65,425.44$  seconds of CPU time to complete the simulation in the sparse-grid computation, while $672,078.25$ seconds of CPU time are needed for finishing the simulation in the corresponding single-grid computation. About $90\%$ CPU time is saved by performing the RK FPFS-WENO simulation on the sparse grids in this 3D problem.

\paragraph{Example 6 (Boat-sail problem).}

In this example, we consider an extension of Voronoi diagram problems as in Example 5,
boat-sail problems (see e.g. \cite{NS2003, NS2005, CD1985}),
which applies a flow field to a Voronoi diagram problem. An application of this kind of problems would be a boat trying to reach the nearest harbor or island on a moving river.

\begin{figure}
	\centering
	\begin{subfigure}[b]{0.45\textwidth}
		\includegraphics[width=\textwidth]{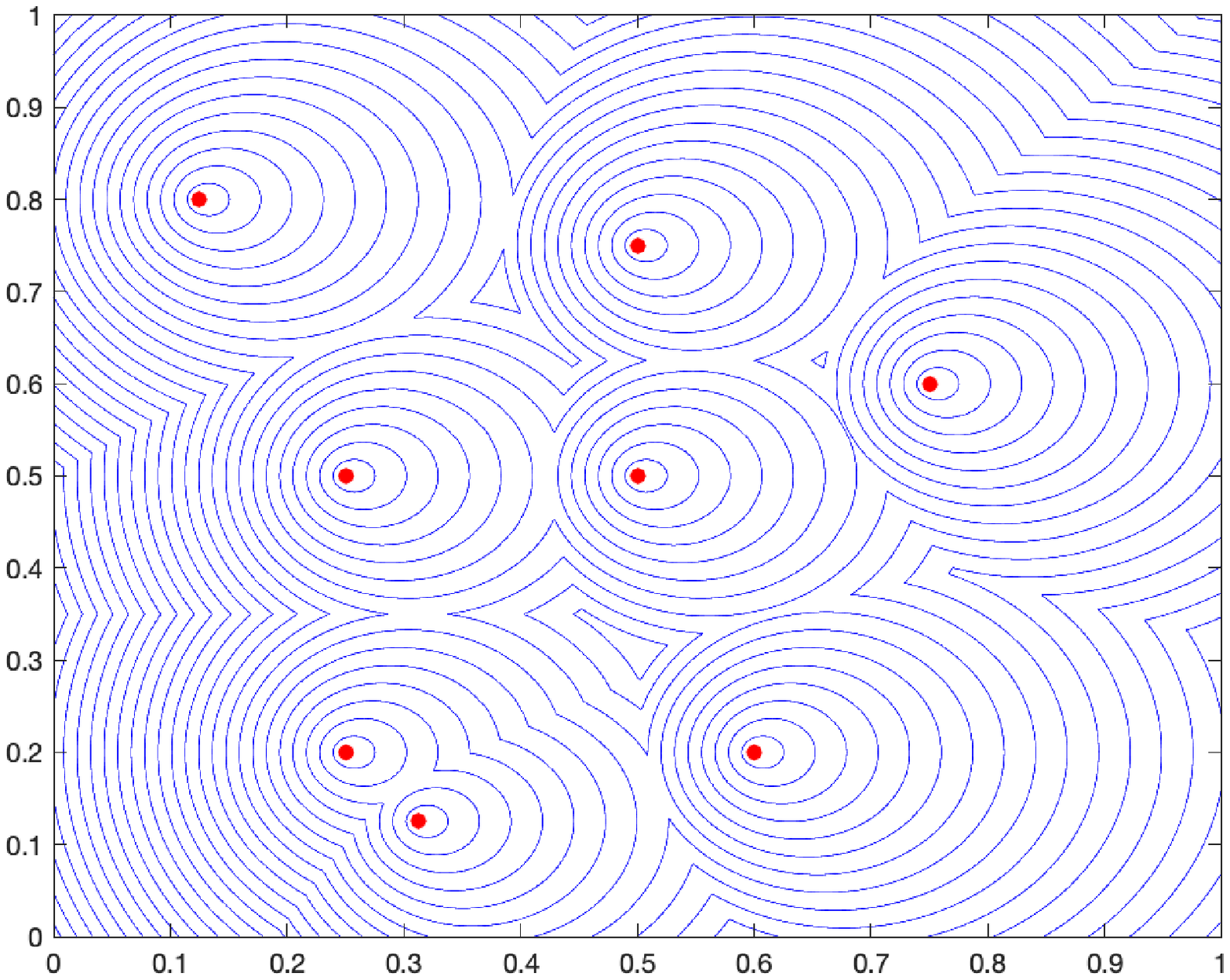}
		\caption{single-grid result}
	\end{subfigure}
	\begin{subfigure}[b]{0.45\textwidth}
		\includegraphics[width=\textwidth]{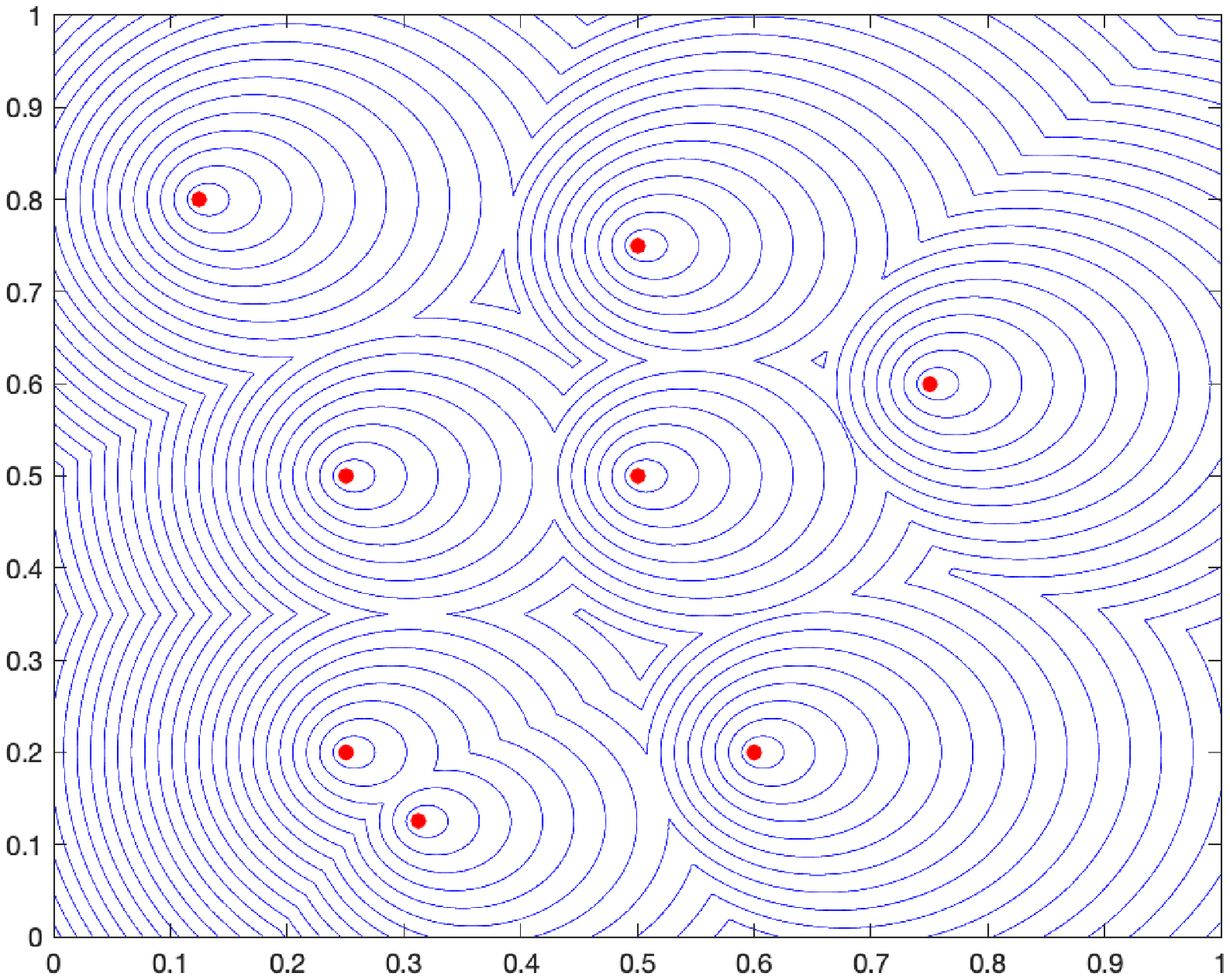}
		\caption{sparse-grid result}
	\end{subfigure}
	\caption{\footnotesize{Example 6, Case 1, numerical solutions of the 2D boat-sail problem by the RK FPFS-WENO scheme on sparse grids ($N_r = 160$ for root grid, finest level $N_L=3$ in the sparse-grid computation) and the corresponding $1280\times1280$ single grid, using the third order WENO interpolation for prolongation in the sparse-grid combination. The contour plots,  30 equally spaced contour lines from $\phi = 0$ to $\phi = 0.4461$. Red points are the harbor locations. (a): single-grid result; (b): sparse-grid result.}}
	\label{fig:2dboat}
\end{figure}

\textbf{Case 1 (2D).} Suppose that the river flows with some velocity $\mathbf{f} = (f_1,f_2)^T$, and the boat travels at a maximum speed $F$ such that $F > |\mathbf{f}|$. The minimum travel time $\phi(x,y)$ from a point $(x,y)$ to the nearest harbor or island can be found by solving the following static Hamilton-Jacobi equation
\begin{align}
F |\nabla \phi| + \mathbf{f} \cdot \nabla \phi &= 1, \quad \mathbf{x} \in \Omega \backslash \Gamma, \label{eq:bs1}\\
\phi(\mathbf{x}) &= 0, \quad \mathbf{x} \in \Gamma,  \label{eq:bs2}
\end{align}
where  $\Gamma$ is the locations of the harbors and islands. Here we take $F = 1$ and $\mathbf{f} = (0.4,0)^T$.
The computational domain is $\Omega = [0,1]^2$ and harbor locations are
$$
\Gamma = \left\{ \left(\frac{1}{4}, \frac{1}{5}\right), \left(\frac{5}{16}, \frac{1}{8}\right), \left(\frac{3}{5}, \frac{1}{5}\right), \left(\frac{3}{4}, \frac{3}{5}\right),  \left(\frac{1}{2}, \frac{3}{4}\right), \left(\frac{1}{4}, \frac{1}{2}\right), \left(\frac{1}{8}, \frac{4}{5}\right), \left(\frac{1}{2}, \frac{1}{2}\right) \right\}.
$$
The exact solution of the problem is also {\em not} smooth.
The sparse grid RK FPFS-WENO scheme with $\gamma = 0.8$ and the third order WENO approximations to the derivatives is applied. The third order WENO interpolation is used for prolongation in sparse-grid computations. Simulations are performed on both sparse grids with $N_r = 160, N_L=3$ and the corresponding $1280\times 1280$ single grid, to compare their numerical results.
The obtained results are reported in Fig.~\ref{fig:2dboat}. As the previous examples, we observe that the numerical solutions by the sparse grid RK FPFS-WENO scheme and its corresponding single-grid computation are comparable, and the nonlinear stability and high resolution properties of the RK FPFS-WENO scheme for resolving the non-smooth solution in this example are preserved well in the sparse-grid simulation. About computational efficiency, it takes $892.71$  seconds of CPU time to complete the simulation in the sparse-grid computation, while $2,155.63$ seconds of CPU time are needed for finishing the simulation in the corresponding single-grid computation. About $59\%$ CPU time is saved by performing the RK FPFS-WENO simulation on the sparse grids in this 2D example.

\begin{figure}
	\centering
	\begin{subfigure}[b]{0.45\textwidth}
		\includegraphics[width=\textwidth]{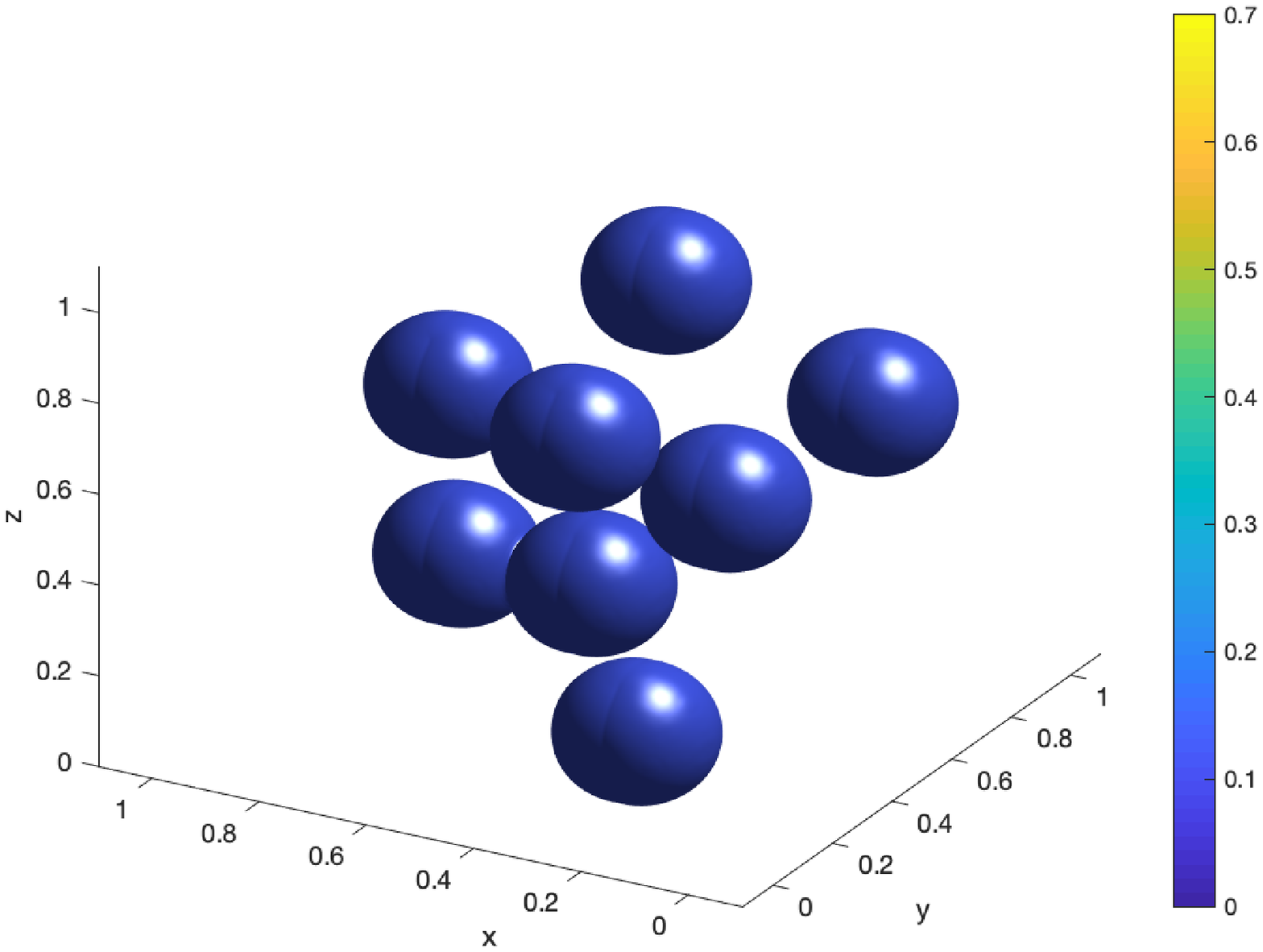}
		\caption{$\phi = 0.125$}
	\end{subfigure}
	\begin{subfigure}[b]{0.45\textwidth}
		\includegraphics[width=\textwidth]{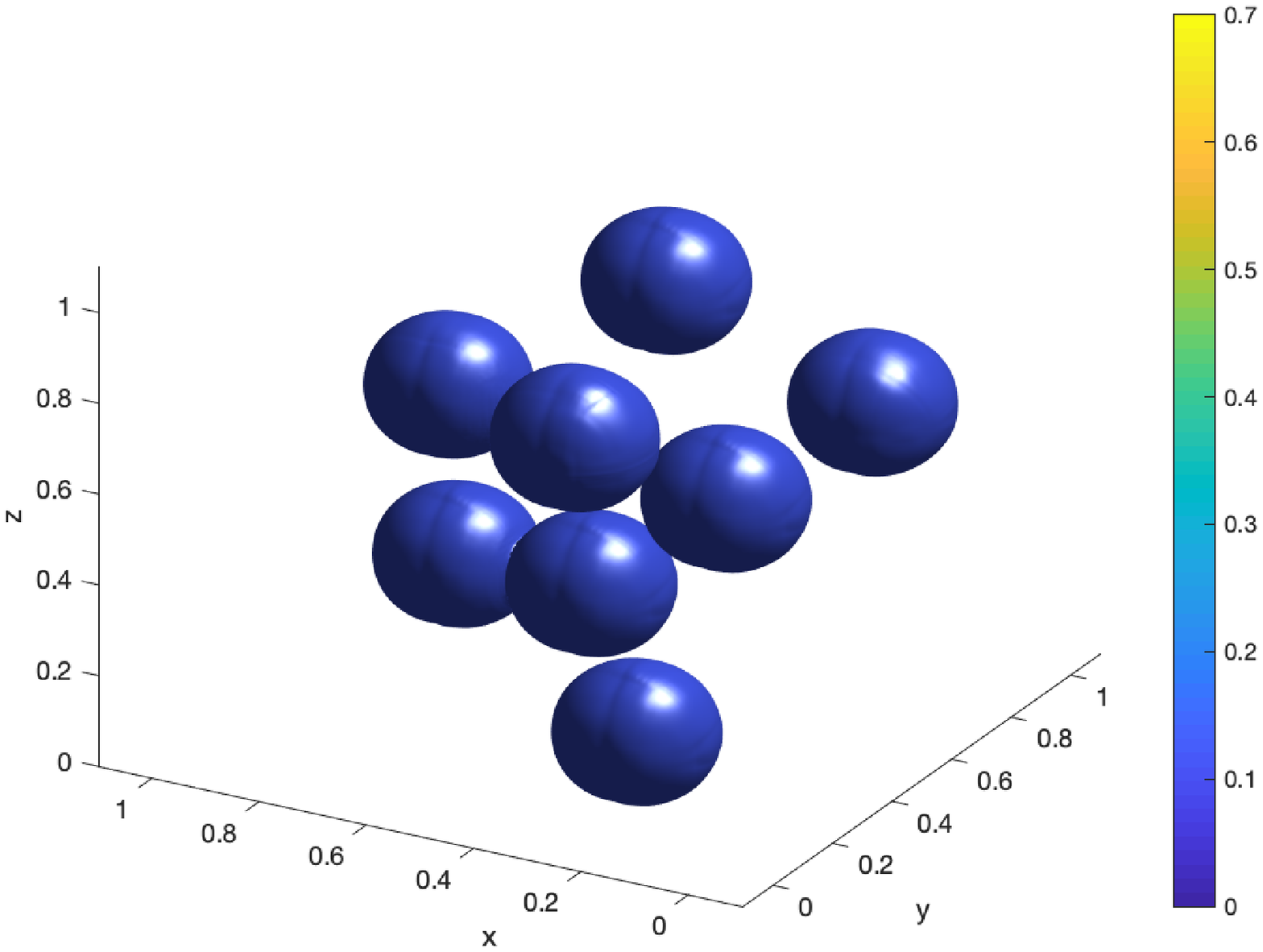}
		\caption{$\phi = 0.125$}
	\end{subfigure}

	\begin{subfigure}[b]{0.45\textwidth}
		\includegraphics[width=\textwidth]{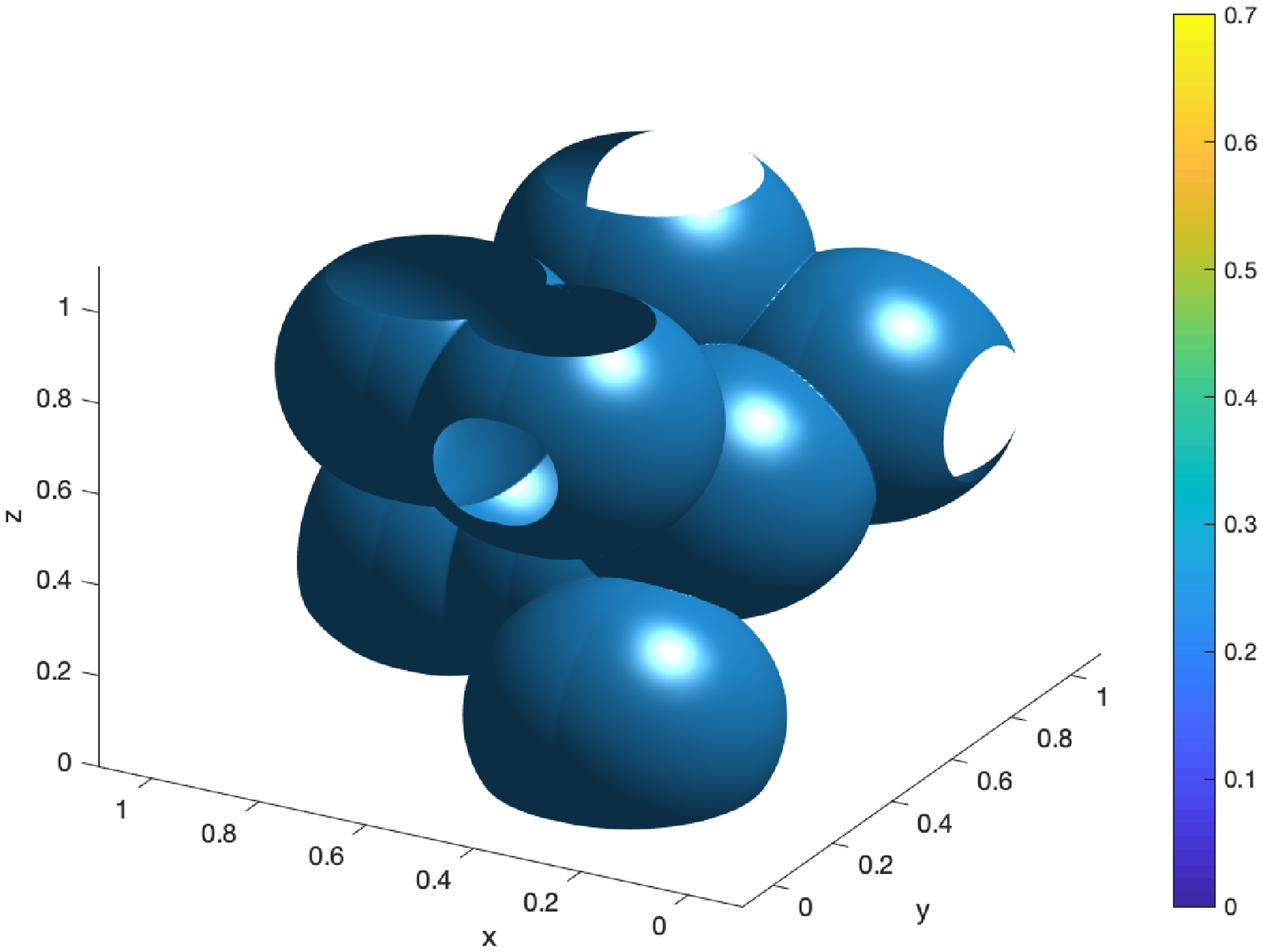}
		\caption{$\phi = 0.25$}
	\end{subfigure}
	\begin{subfigure}[b]{0.45\textwidth}
		\includegraphics[width=\textwidth]{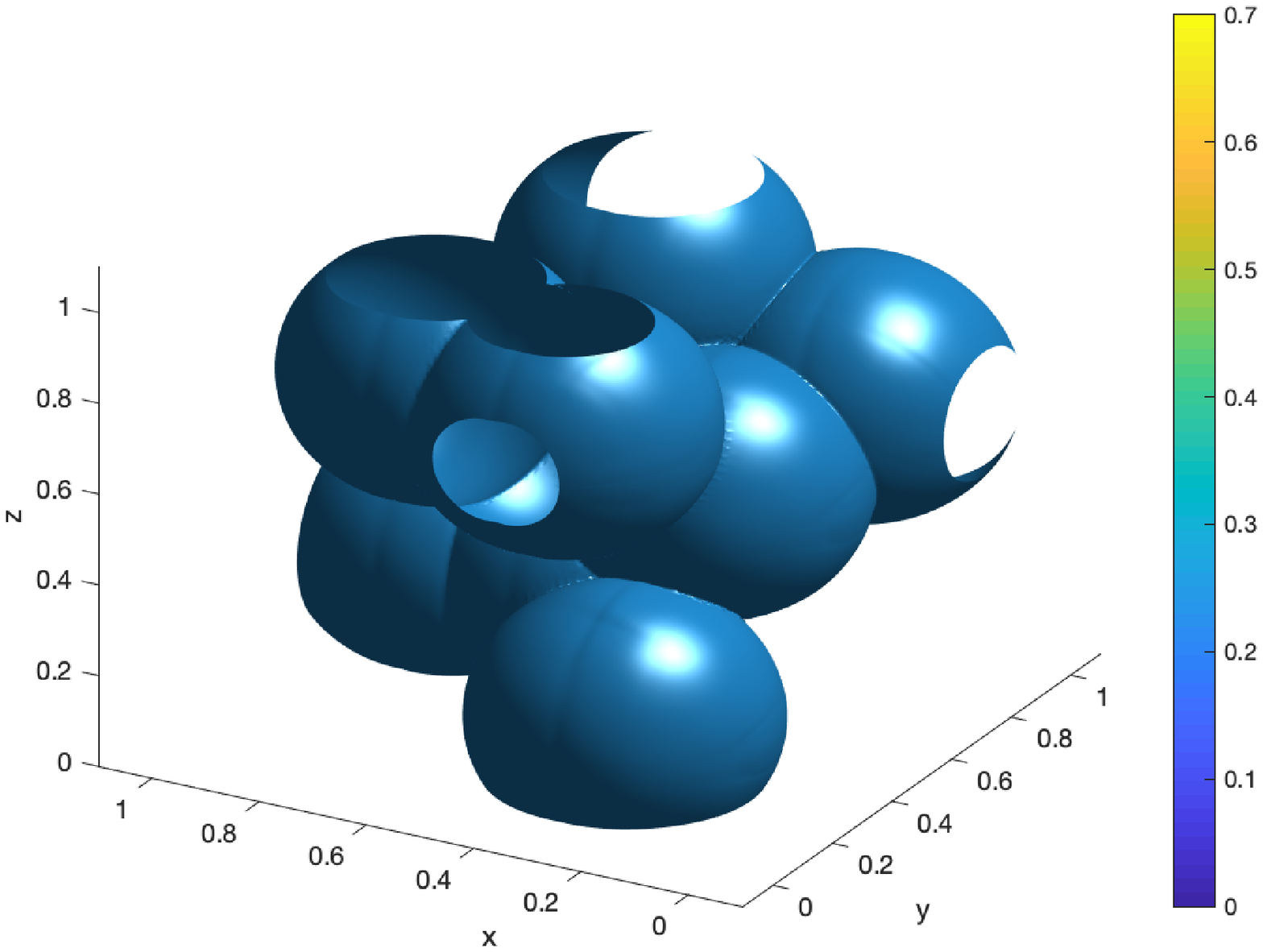}
		\caption{$\phi = 0.25$}
	\end{subfigure}

	\begin{subfigure}[b]{0.45\textwidth}
		\includegraphics[width=\textwidth]{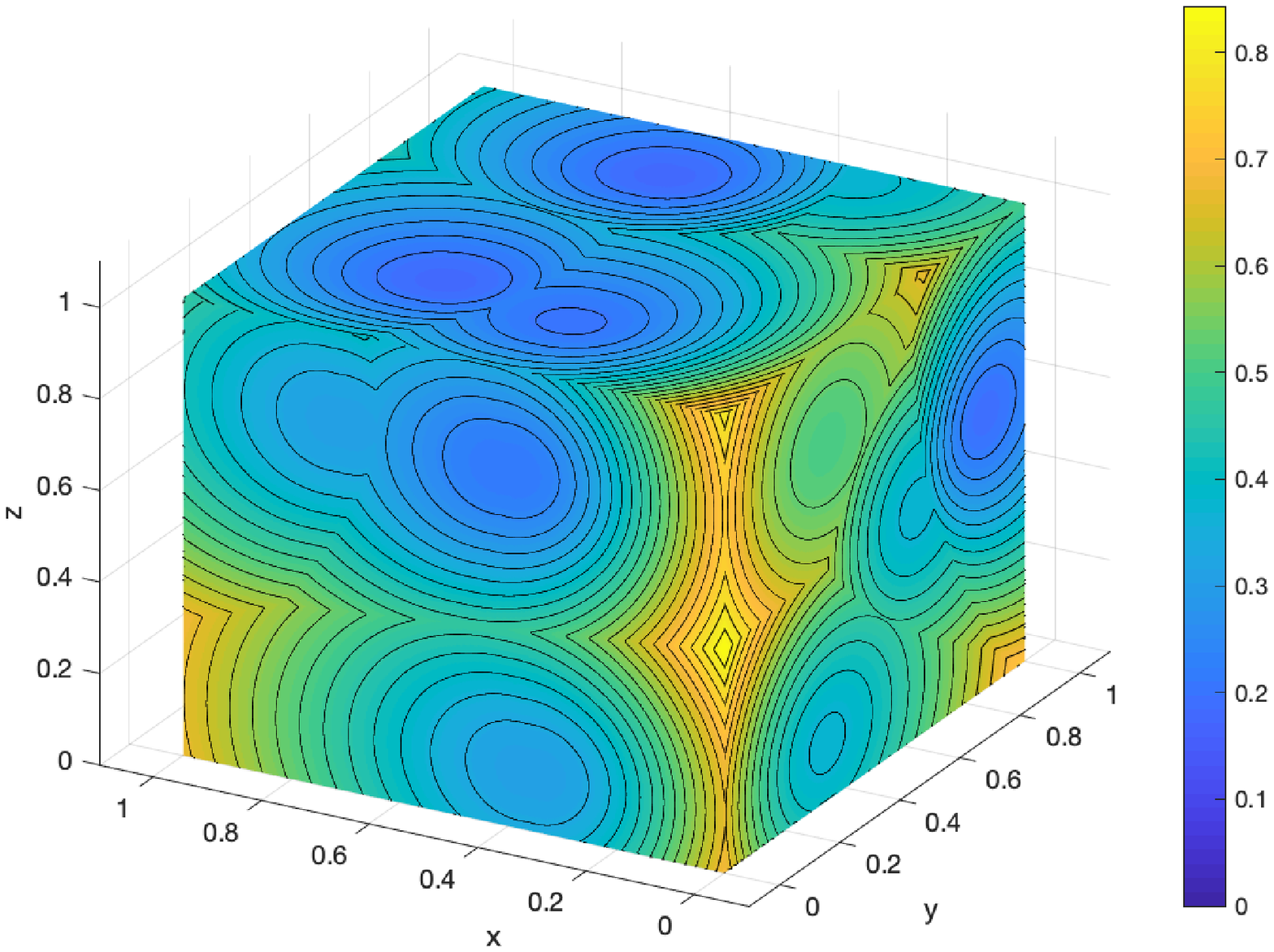}
		\caption{Surface contour}
	\end{subfigure}
	\begin{subfigure}[b]{0.45\textwidth}
		\includegraphics[width=\textwidth]{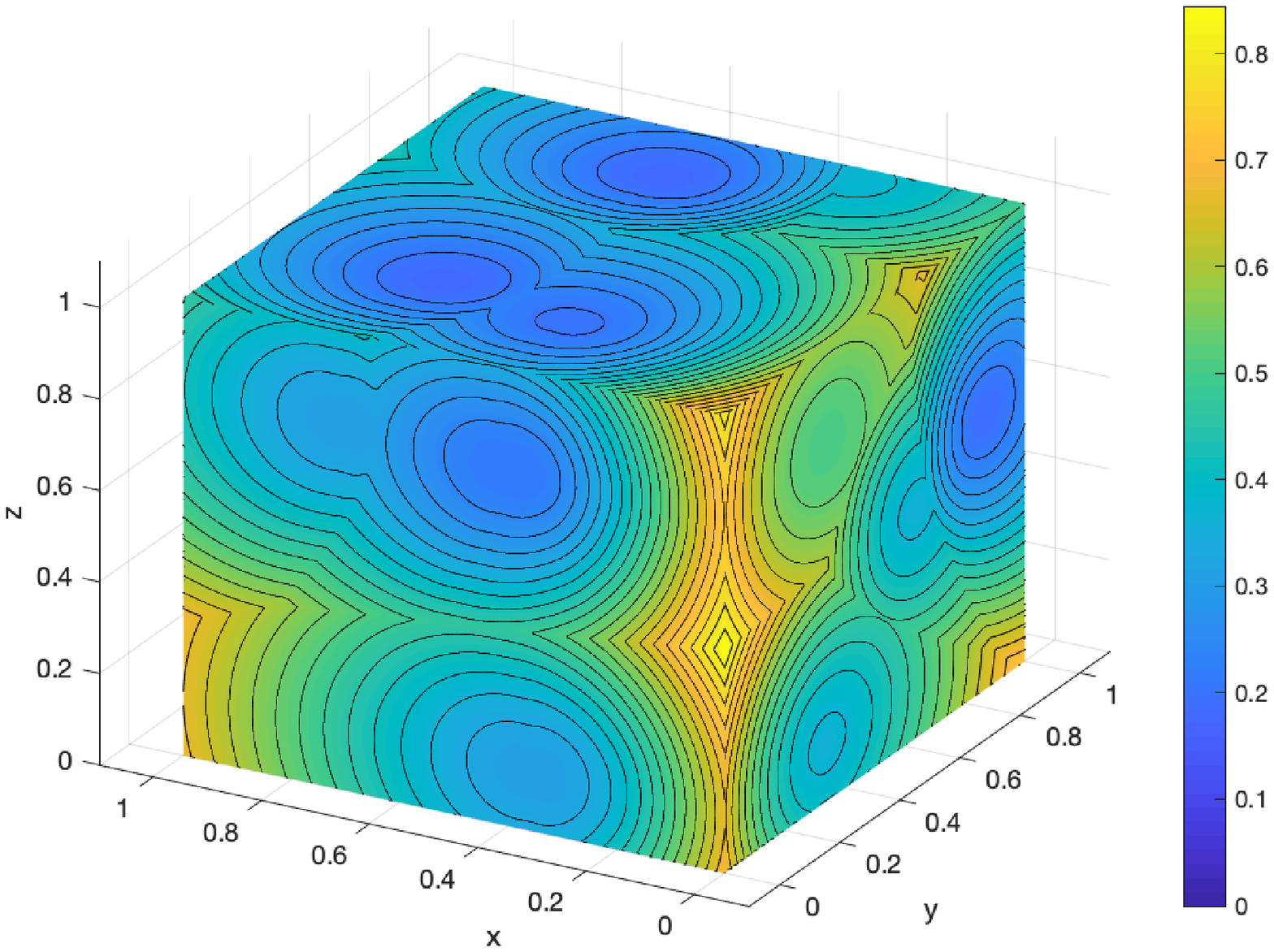}
		\caption{Surface contour}
	\end{subfigure}
	\caption{\footnotesize{Example 6, Case 2, numerical solutions of the 3D boat-sail problem by the RK FPFS-WENO scheme on sparse grids ($N_r = 80$ for root grid, finest level $N_L=3$ in the sparse-grid computation) and the corresponding $640\times640\times640$ single grid, using the third order WENO interpolation for prolongation in the sparse-grid combination. (a), (c), (e): single-grid result; (b), (d), (f): sparse-grid result; (a), (b): the contour plots for $\phi=0.125$; (c), (d): the contour plots for $\phi=0.25$; (e), (f): the contour plots for the whole surface.}}
	\label{fig:3dboat}
\end{figure}

\textbf{Case 2 (3D).} We consider the 3D case of the problem. The river flows with velocity $\mathbf{f} = (f_1,f_2,f_3)^T$, and the boat travels at a maximum speed $F$ such that $F > |\mathbf{f}|$. The minimum travel time $\phi(x,y,z)$ from a point $(x,y,z)$ to the nearest harbor or island can be found by solving the 3D version of the static Hamilton-Jacobi equation~\eqref{eq:bs1}-\eqref{eq:bs2}. We take $F = 1$ and $\mathbf{f} = (0.4,0.4,0)^T$.
The computational domain is $\Omega = [0,1]^3$ and harbor locations are
\begin{align*}
\Gamma &= \left\{ \left( \frac{1}{4}, \frac{1}{5}, \frac{1}{8}\right), \left( \frac{1}{3}, \frac{1}{7}, \frac{7}{9}\right), \left( \frac{3}{5}, \frac{1}{5}, \frac{4}{5}\right), \left( \frac{3}{4}, \frac{1}{2}, \frac{1}{4}\right), \right. \\
&\qquad \left. \left( \frac{1}{2}, \frac{3}{4}, \frac{4}{5}\right), \left( \frac{1}{4}, \frac{1}{2}, \frac{1}{2}\right), \left( \frac{1}{7}, \frac{4}{5}, \frac{3}{5}\right), \left( \frac{1}{2}, \frac{1}{2}, \frac{1}{4}\right) \right\}.
\end{align*}
Again, the exact solution of the 3D problem is {\em not} smooth. The sparse grid RK FPFS-WENO scheme with $\gamma = 0.8$ and the third order WENO approximations to the derivatives is applied in solving this 3D problem. The third order WENO interpolation is used for prolongation in sparse-grid computations. Simulations are carried out on both sparse grids with $N_r = 80, N_L=3$ and the corresponding $640\times 640\times 640$ single grid, for comparing the numerical results.
The simulation results are presented in Fig.~\ref{fig:3dboat}, which shows that the numerical solutions by the sparse grid RK FPFS-WENO scheme and its corresponding single-grid computation are comparable, and the nonlinear stability and high resolution properties of the RK FPFS-WENO scheme for resolving the non-smooth solution in this 3D boat-sail problem are preserved well in the sparse-grid simulation. In terms of computational efficiency, it takes $179,925.93$  seconds of CPU time to complete the simulation in the sparse-grid computation, while $1,480,989.58$ seconds of CPU time are needed for finishing the simulation in the corresponding single-grid computation. About $88\%$ CPU time is saved by performing the RK FPFS-WENO simulation on the sparse grids in this 3D boat-sail problem.

\section{Conclusions}

In this technical note, we apply the sparse-grid combination technique to a third order  Runge-Kutta type fixed-point fast sweeping WENO (RK FPFS-WENO) scheme for efficiently computing solutions of multidimensional Eikonal equations.
Due to their sophisticated nonlinearity, more computational costs than
many other schemes are needed in high order WENO simulations, especially for multidimensional problems. Here we follow our previous work and implement the RK FPFS-WENO scheme on sparse grids. A third order WENO interpolation is applied in the prolongation step of the sparse-grid combination technique, for robust computations of non-smooth solutions of Eikonal equations in sparse-grid simulations. Numerical experiments on 2D and 3D problems are performed for the sparse grid RK FPFS-WENO method to show that a more efficient algorithm than regular RK FPFS-WENO method on single grids to solve the multidimensional Eikonal equations is achieved, with about $50\% \sim 90\%$  CPU time costs being saved on refined meshes, by comparing with the corresponding single-grid simulations in examples here.

In this technical note, we focus on the efficient implementation of the RK FPFS-WENO scheme on sparse grids and its numerical experiments. We would like to point out that there are still quite a few open problems to be investigated further for the sparse grid method. For example, it is still an open problem on how to perform theoretical error analysis for such kind of nonlinear sparse grid schemes, although that has been done for the linear schemes in solving linear time-dependent PDEs in the literature. In the numerical experiments, we find that
the $L^\infty$ errors of sparse-grid computations are larger than the corresponding single-grid computations.
How to improve the accuracy of the sparse grid scheme is an interesting and important problem. Furthermore, the
sparse grid RK FPFS-WENO method studied here should be able to be extended to higher order accuracy and to solve more complicated static HJ equations. All of these important open problems
will be our future work.

\bigskip
\bigskip
\bigskip
\noindent {\bf Conflict of Interest:} the authors declare that there is no conflict of interest.

\end{document}